\documentclass[12pt,reqno]{amsart}
\usepackage{amsmath,amsfonts,amsthm,amssymb}
\newcommand\version{October 26, 2006}
\usepackage{amsxtra}

\font\notefont=cmsl8  
\newtheorem{theorem}{Theorem}[section]
\newtheorem{proposition}[theorem]{Proposition}
\newtheorem{lemma}[theorem]{Lemma}
\newtheorem{corollary}[theorem]{Corollary}

\theoremstyle{definition}

\theoremstyle{remark}
\newtheorem{remark}[theorem]{Remark}

\numberwithin{equation}{section}

\newcommand{\C}{\mathbb{C}}

\newcommand{\eps}{\epsilon}
\renewcommand{\epsilon}{\varepsilon}

\newcommand{\N}{\mathbb{N}}

\renewcommand{\phi}{\varphi}
\newcommand{\R}{\mathbb{R}}
\newcommand{\Sph}{\mathbb{S}}

\DeclareMathOperator{\dom}{dom}

\DeclareMathOperator{\re}{Re}
\DeclareMathOperator{\supp}{supp}

\DeclareMathOperator{\tr}{tr}

\def\slim{\mathop{\mbox{s-lim}}\limits}

\begin{document}

\title[Hardy-Lieb-Thirring inequalities -- \version]{Hardy-Lieb-Thirring
     inequalities for fractional Schr\"odinger operators}

\author[Rupert L. Frank]{Rupert L. Frank}
\address{Rupert L. Frank, Department of Mathematics, Royal Institute
     of Technology, 100 44 Stockholm, Sweden}
\email{rupert@math.kth.se}

\author[Elliott H. Lieb]{Elliott H. Lieb}
\address{Elliott H. Lieb, Departments of Mathematics and Physics,
Princeton University,
     P.~O.~Box 708, Princeton, NJ 08544, USA}
\email{lieb@princeton.edu}

\author[Robert Seiringer]{Robert Seiringer}
\address{Robert Seiringer, Department of Physics, Princeton University,
P.~O.~Box 708,
     Princeton, NJ 08544, USA}
\email{rseiring@princeton.edu}

\begin{abstract}
     We show that the Lieb-Thirring inequalities on moments of negative
     eigenvalues of Schr\"odinger-like operators remain true, with
     possibly different constants, when the critical Hardy-weight $C
     |x|^{-2}$ is subtracted from the Laplace operator. We do so by first
     establishing a Sobolev inequality for such operators. Similar
     results are true for fractional powers of the Laplacian and the
     Hardy-weight and, in particular, for relativistic Schr\"odinger
     operators. We also allow for the inclusion of magnetic vector
     potentials. As an application, we extend, for the first time, the
     proof of stability of relativistic matter with magnetic fields all
     the way up to the critical value of the nuclear charge
     $Z\alpha=2/\pi$.
\end{abstract}

\thanks{\copyright\, 2006 by the authors. This paper may be reproduced, in its
entirety, for non-commercial purposes.}

\date{\version}
\maketitle

\vspace{-.7cm}
{\footnotesize{\tableofcontents}}

\section{Introduction}

We shall generalize several well known
inequalities about the negative spectrum of Schr\"odinger-like
operators on $\R^d$. As an application of our results we shall give a
proof of the `stability of relativistic matter' --- one which goes
further than previous proofs by permitting the inclusion of magnetic
fields for values of the nuclear charge all the way up to $Z\alpha =
2/\pi$, which is the critical value in the absence of a field.

There are three main inequalities to which we refer.  The
first is Hardy's inequality, whose classical form for $d\geq 3$ is the
following. (In this introduction we shall not be precise about the
space of functions in question, but will be precise later on.)
\begin{equation} \label{hardy} \int_{\R^d} \left (|\nabla u(x)|^2 -
\frac{(d-2)^2}{4|x|^2} |u(x)|^2 \right)
    dx  \geq 0 \ .
\end{equation}

The second is Sobolev's inequality for $d\geq 3$, \begin{equation} \label{sob}
    \int_{\R^d} |\nabla u(x)|^2 \, dx\geq S_{2,d}\ \left\{ 
\int_{\R^d}|u(x)|^{2d/(d-2)} dx \right\}^{1-2/d}
=  S_{2,d} \Vert u \Vert_{2d/(d-2)}^2    \ .
\end{equation}

The third is the Lieb-Thirring (LT) inequality \cite{LTh} for the Schr\"odinger
operator $H=-\Delta - V(x)$. If its negative eigenvalues are denoted
by $-\lambda_1\leq -\lambda_2\leq \cdots$, and if $\gamma \geq 0$,
then
\begin{equation} \label{lt}
    \sum_{j}\lambda_j^\gamma \leq L_{\gamma,d}\int_{\R^d}
V(x)_+^{\gamma+d/2}\,dx
= L_{\gamma,d} \Vert V_+\Vert_{ \gamma+d/2}^{   \gamma+d/2}\,  .
\end{equation}
This holds if and only if $\gamma\geq\frac12$ when $d=1$, $\gamma>0$ when
$d=2$ and $\gamma\geq 0$ when $d\geq3$. (Here and in the sequel
$t_-:=\max\{ 0, -t\}$ and $t_+:=\max\{ 0,t\}$ denote the negative and
positive parts of $t$.)

By duality, \eqref{sob} is equivalent to the fact that the
Schr\"odinger operator $H = -\Delta -V$, with has no
negative eigenvalues if $\Vert V_+\Vert_{d/2} \leq S_{2,d}$. On the other
hand, \eqref{lt} gives an upper bound to the number of negative
eigenvalues in terms of $\Vert V_+\Vert_{d/2}$ when $\gamma=0$ and it
estimates the magnitude of these eigenvalues when $\gamma >0$.

All three inequalities can be generalized by the inclusion of a magnetic vector
potential $A$ (related to the magnetic field $B$ by $B={\rm curl}A$).
That is, $\nabla $ is replaced by $\nabla - i A(x)$, and $\Delta$ by
$(\nabla - i A(x))^2$. The sharp constants in \eqref{hardy},
\eqref{sob} remain unchanged while the sharp constants in \eqref{lt}
that are independent of $A$ might, in principle, be different from the ones for
$A=0$. However, the best constants known so far do not depend on $A$.
The inclusion of $A$ is easily done in \eqref{hardy}, \eqref{sob} by
using the diamagnetic inequality, but the inclusion in \eqref{lt} is
more delicate; one uses the Feynman-Kac path integral formula to show
that for each $x,y \in \R^d$ and $t, \tau>0$, the $A$-field reduces
the magnitude of the heat kernel $e^{t(\nabla - i A)^2}(x,y)$ relative
to $e^{t\Delta}(x,y)$, and hence reduces the resolvent kernel $ |[-(\nabla -
i A)^2 + \tau ]^{-1}(x,y)|$ relative to $[-\Delta +\tau ]^{-1}(x,y)$.

{}From \cite{BVa} one can deduce that  \eqref{hardy} and \eqref{sob} can
be combined as follows: For $\eps>0$, \begin{equation} \label{frank1}
    \int_{\R^d} \left(|\nabla u(x)|^2 -  \frac{(d-2)^2}{4|x|^2} |u(x)|^2
+ |u(x)|^2  \right) dx \geq S'_{2,d,\epsilon}\ \Vert u \Vert_{2d/(d-2 
+\epsilon)}^2 \, .
\end{equation}
with $S'_{2,d,\epsilon} \to 0$ as $\epsilon \to 0$.  Note the
extra term $|u(x)|^2$ on the left side to account for the fact that
the left and right sides behave differently under scaling; examples show
that it really is necessary to have $\epsilon>0$ here.

In \cite{EkFr}, a parallel extension of \eqref{lt} for the negative
eigenvalues $-\lambda_j$ of the Schr\"odinger
operator $H = -\Delta - (d-2)^2/(4|x|^2)- V$ is proved. For $\gamma>0$ and
$d\geq3$, \begin{equation} \label{lt-hardy}
    \sum_{j}\lambda_j^\gamma \leq L'_{\gamma,d}\int_{\R^d}
V(x)_+^{\gamma+d/2}\,dx
=L'_{\gamma,d}\Vert V_+\Vert_{ \gamma+d/2}^{   \gamma+d/2}\ ,             \,
\end{equation}
with $L'_{\gamma,d} \geq L_{\gamma,d}$. Note that there is no need
for an $\epsilon$ in (\ref{lt-hardy}); the fact that $\epsilon \neq 0$ in 
\eqref{frank1}
is reflected here in the fact that $\gamma >0$ is needed. As before,
a magnetic vector potential can easily be included in \eqref{frank1},
but it does not seem easy to include a magnetic field in \eqref{lt-hardy}
by the methods in \cite{EkFr}.

Our goal here is to extend these results in several ways. One
extension is to include a magnetic field in \eqref{lt-hardy}. Another
is to consider fractional powers of the (magnetic) Laplacian, i.e., to
the case in which $|\nabla - i A|^2$ is replaced by $|\nabla - i
A|^{2s}$ with $0< s <\min\{1,d/2\}$ (which means that we can now
include one- and two-dimensions).  This is a significant
generalization because the operator $(-\Delta)^{s}$ is not a
differential operator and it is not `local'. Really different
techniques will be needed.  In particular, we shall use the heat
kernel to prove the analog of \eqref{lt-hardy}, in the
manner of \cite{L1}. A bound on this kernel, in turn, will be derived from
a  Sobolev-like inequality (the analogue of \eqref{frank1}) by using
an analogue of Nash's inequality, as explained in \cite{LLo}. The
appropriate inequalities are naturally formulated in a weighted space
with measure $|x|^{-\beta} dx$ for $\beta>0$. Therefore, a pointwise
bound on the heat kernel for a weighted `Hardy' operator
$\exp\{-t |x|^\alpha ((-\Delta)^s - \mathcal
C_{s,d}|x|^{-2s}+1)|x|^\alpha \}$ for appropriate $\alpha>0$ will be
needed and will not be straightforward to obtain.

In the dimension most relevant for physics, $d=3$, the earlier case
$s=1$ may be called the {\it non-relativistic} case, while the new
result for $s=1/2$ may be called the {\it relativistic} case. Indeed,
the resulting LT inequality, together with some of the methodology in
\cite{LY}, yields a new proof of the stability of relativistic matter,
which will be sketched in Section~\ref{sec:stability}. The main point,
however, is that this new proof allows for an arbitrary magnetic
vector potential $A$.  Since the constant in the relativistic ($s=1/2$) Hardy
inequality that replaces $(d-2)^2/4$ is  $2/\pi$ (which is the
same as the critical value of $Z\alpha$ in the field-free relativistic
case), we conclude that we can simultaneously have an arbitrary $A$-field
and the critical value of $Z\alpha$, the nuclear charge times the
fine-structure constant.  Up to now it was not possible to have both
an arbitrary field and $Z\alpha =2/\pi$. Therefore, the proof of the
analogue of \eqref{lt-hardy} for $s=1/2$ with an $A$-field opens a
slightly improved perspective on the interaction of matter and
radiation.

\bigskip {\it Acknowledgements.} We thank Heinz Siedentop for
suggesting that we study inequalities of this type, and we thank him,
Ari Laptev and Jan Philip Solovej for helpful discussions. We also
thank Renming Song for valuable comments on a previous version of
this manuscript. This work was partially supported by the Swedish
Foundation for International Cooperation in Research and Higher
Education (STINT) (R.F.), by U.S.  National Science Foundation grants
PHY 01 39984 (E.L.) and PHY 03 53181 (R.S.), and by an A.P. Sloan
Fellowship (R.S.).


\section{Main Results}

\subsection{Hardy-Lieb-Thirring inequalities}

We recall the Hardy-type inequality
\begin{equation}\label{eq:hardy}
     \int_{\R^d} |x|^{-2s} |u(x)|^2 \,dx
     \leq \mathcal C_{s,d}^{-1}
     \int_{\R^d} |\xi|^{2s} |\hat u(\xi)|^2 \,d\xi,
     \qquad u\in C_0^\infty(\R^d),
\end{equation}
valid for $0<2s<d$. Here
\begin{equation*}
     \hat u(\xi) := (2\pi)^{-d/2} \int_{\R^d} u(x) e^{-i\xi\cdot x}\,dx,
     \qquad \xi\in\R^d,
\end{equation*}
denotes the Fourier transform of $u$. The sharp constant in
\eqref{eq:hardy},
\begin{equation} \label{eq:csd}
     \mathcal C_{s,d} := 2^{2s} \frac{\Gamma^2((d+2s)/4)}{\Gamma^2((d-2s)/4)}\ ,
\end{equation}
has been found independently by Herbst \cite{He} and Yafaev
\cite{Ya}. Moreover, in the latter paper it is shown that this
constant is not achieved in the class of functions for which both
sides are finite. In the case $0<s<\min\{1,d/2\}$, this fact can also
be deduced from our  ground state representation in Proposition
\ref{substitution}, which therefore represents an independent proof of
\eqref{eq:hardy}.   See Remark~\ref{rem:hardy} below.

We denote $D=-i\nabla$. Consider a magnetic vector potential $A\in
L^2_{\rm loc}(\R^d,\R^d)$ and the self-adjoint operator $(D-A)^2$ in
$L^2(\R^d)$. For $0<s\leq 1$ we define the operator
$|D-A|^{2s}:=((D-A)^2)^s$ by the spectral theorem. One form of the
diamagnetic inequality states that if $0<s\leq 1$ and $u\in
\mathcal\dom |D-A|^{s}$ then $|u|\in H^s(\R^d)$ and
\begin{equation}\label{eq:diamag1}
           \big\|(-\Delta)^{s/2} |u|\big\|^2
     \leq \big\||D-A|^{s} u\big\|^2.
\end{equation}
Here and in the sequel, $\|\cdot\|=\|\cdot\|_2$ denotes the $L^2$-norm. We
refer to Remark \ref{diamag1} for more details concerning
\eqref{eq:diamag1}. Combining this inequality with
the Hardy inequality \eqref{eq:hardy} we find that the quadratic form
\begin{equation}\label{eq:defform}
     h_{s,A}[u] := \big\| |D-A|^{s} u\big\|^2- \mathcal C_{s,d} \big\| 
|x|^{-s}u\big\|^2
\end{equation}
is non-negative on $\dom |D-A|^{s}$ if $0<s<\min\{1,d/2\}$. We use the
same notation for its closure and denote by
\begin{equation*}
     H_{s,A} = |D-A|^{2s} - \mathcal C_{s,d}|x|^{-2s}
\end{equation*}
the corresponding self-adjoint operator in $L^2(\R^d)$.

Our main result is

\begin{theorem}[{\bf Hardy-Lieb-Thirring inequalities}]\label{main}
     Let $\gamma>0$ and $0<s<\min\{1,d/2\}$. Then
     there exists a constant $L_{\gamma,d,s}>0$ such that for all $V$ and $A$
     \begin{equation}\label{eq:main}
       \tr\left(|D-A|^{2s} - \mathcal C_{s,d}|x|^{-2s} -
       V\right)_-^\gamma
       \leq L_{\gamma,d,s} \int_{\R^d} V(x)_+^{\gamma+d/2s} \,dx.
     \end{equation}
\end{theorem}

\begin{remark}\label{rem:s1}
     For $d\geq 3$, Theorem~\ref{main} holds also for $s=1$. In the
     non-magnetic case $A=0$, this has been proved in \cite{EkFr}. In the
     present paper, we present an independent proof of this result, which
     allows for the inclusion of an $A$-field. For simplicity, we
     restrict our attention to $s<1$ in the following, and comment on the
     (simpler) case $s=1$ in Subsection~\ref{sec:s1}.
\end{remark}

Theorem~\ref{main} will be proved in Section~\ref{sec:lth}.
The main ingredient in its proof is a Sobolev-type
inequality which might be of independent interest and which we shall
present in the remainder of this subsection.

For the case $A=0$ we shall drop the index $A$ from the notation, i.e.,
\begin{equation*}
     h_s[u] =
     \int_{\R^d} |\xi|^{2s} |\hat u(\xi)|^2 \,d\xi
     - \mathcal C_{s,d} \int_{\R^d} |x|^{-2s} |u(x)|^2 \,dx.
\end{equation*}
We denote the closure of this form by the same letter. In particular,
its domain $\dom h_s$ is the closure of $C_0^\infty(\R^d)$ with
respect to the norm $(h_s[u]+\|u\|^2)^{1/2}$. Note that
$H^s(\R^d)\subset\dom h_s$ with \emph{strict} inclusion. In
particular, there exist functions $u\in\dom h_s$ for which both sides
of \eqref{eq:hardy} are infinite.

Hardy's inequality~(\ref{eq:hardy}) implies that $h_s$ is
non-negative. The following theorem shows that for functions of compact
support it even satisfies a Sobolev-type inequality; i.e., $h_s[u]$
can be bounded from below by an $L^q$-norm of $u$.

\begin{theorem}[{\bf Local Sobolev-Hardy inequality}]\label{sobolev}
     Let $0<s<\min\{1,d/2\}$ and $1\leq q <
     2^* = 2d/(d-2s)$. Then there exists a constant $C_{q,d,s}>0$ such that
for any domain $\Omega\subset\R^d$ with finite measure $|\Omega|$ one has
     \begin{equation}\label{eq:sobolev}
       \|u\|_q^2 \leq C_{q,d,s} |\Omega|^{2\left(\frac1q-\frac1{2^*}\right)}
h_s[u],
       \qquad u\in C_0^\infty(\Omega).
     \end{equation}
\end{theorem}

\begin{remark}
     Note that the exponent $q$ is \emph{strictly smaller} than the
     critical Sobolev exponent $2^*=2d/(d-2s)$. By considering functions
     which diverge like $|x|^{-(d-2s)/2}$ at $x=0$, it is easy to see that
     the inequality \eqref{eq:sobolev} cannot hold if $q$ is replaced by
     $2^*$.
\end{remark}

We note that the analogue of Theorem \ref{sobolev} in the local case $s=1$
   is proved in \cite{BVa}. Theorem~\ref{sobolev}  will be proved in
Section~\ref{sec:sobolev}, where we also deduce the following corollary from 
it.

\begin{corollary}[{\bf Global Sobolev-Hardy inequality}]\label{sobolevcor}
     Let $0<s<\min\{1,d/2\}$ and $2\leq q < 2^*=
     2d/(d-2s)$. Then there exists a constant $C_{q,d,s}'>0$ such that
     \begin{equation}\label{eq:sobolevcor}
       \|u\|_q^2 \leq C_{q,d,s}' (h_s[u] + \|u\|^2),
       \qquad u\in C_0^\infty(\R^d).
     \end{equation}
\end{corollary}

Note that \eqref{eq:sobolevcor} may be written in the scale-invariant form
\begin{equation}\label{eq:corscale}
           \|u\|_q \leq C_{q,d,s}'' h_s[u]^{\frac 
d{2s}\left(\frac12-\frac1q\right)}
           \|u\|^{\frac ds\left(\frac1q-\frac1{2^*}\right)}\,,
\end{equation}
where $C_{q,d,s}''$ can be expressed explicitly in terms of
$C_{q,d,s}'$. This inequality follows from applying
(\ref{eq:sobolevcor}) to functions of the form $u_\lambda(x)=
u(\lambda x)$ and then optimizing over the choice of $\lambda$.

Although the sharp constants in Theorems~\ref{main} and~\ref{sobolev},
as well as in Corollary~\ref{sobolevcor}, are unknown, explicit upper
bounds involving a certain variational expression can be deduced from
our proof.  In Appendix~\ref{app:sobolevconst}, we do evaluate explicit
bounds on the constants in Theorem~\ref{sobolev} in the special case
$d=3$, $s=1/2$, which is the most interesting case from a physical
point of view, as will be explained in the next subsection.

\begin{remark}
     Corollary~\ref{sobolevcor} is one of the main ingredients in our
     proof of Theorem~\ref{main}. On the other hand,
     Corollary~\ref{sobolevcor} is an easy consequence of
     Theorem~\ref{main}, except for the value of the constants. In fact,
     (\ref{eq:main}) implies that $$
h_{s}[u] \geq \int_{\R^d} V(x) |u(x)|^2\,dx - \left( L_{\gamma,d,s}
     \int_{\R^d} V(x)_+^{\gamma+d/2s} \, dx\right)^{1/\gamma}\|u\|^2
$$
for all $u\in \dom h_s$ and all $V\in
L^{\gamma+d/2s}(\R^d)$. Eq.~(\ref{eq:corscale}) follows by optimizing
the right side over all $V$. \end{remark}


\subsection{Stability of relativistic matter}\label{sec:stability}

We shall now explain how the inequalities in Theorem~\ref{main} can be
used to prove stability of relativistic matter in the presence of an
external magnetic field. The proof works up to and including the
critical value $Z\alpha = 2/\pi$, which is a new result and solves a
problem that has been open for a long time. We refer to \cite{L2,L3}
for a review of this topic.

We consider $N$ electrons of mass $m\geq 0$ with $q$ spin states
($q=2$ for real electrons) and $K$ fixed nuclei with (distinct)
coordinates $R_1,\ldots,R_K\in\R^3$ and charges $Z_1,\ldots,Z_K> 0$. A
pseudo-relativistic description of the corresponding
quantum-mechanical system is given by the Hamiltonian
\begin{equation*}
     H_{N,K} = \sum_{j=1}^N 
\left(\sqrt{(D_j-\sqrt{\alpha}A(x_j))^2+m^2}-m\right)
     + \alpha V_{N,K}(x_1,\dots,x_N)\,.
\end{equation*}
The Pauli exclusion principle dictates that $H_{N,K}$ acts on
functions in the anti-symmetric $N$-fold tensor product of
$L^2(\R^3;\C^q)$. Here we use units where $\hbar=c=1$, $\alpha>0$ is the
fine structure constant, and
\begin{align*}
     V_{N,K}(x_1,\dots,x_N) := & \sum_{1\leq i<j\leq N} |x_i-x_j|^{-1}
     - \sum_{j=1}^N \sum_{k=1}^K Z_k |x_j-R_k|^{-1} \\
     & + \sum_{1\leq k<l\leq K} Z_k Z_l |R_k-R_l|^{-1}.
\end{align*}
Stability of matter means that $H_{N,K}$ is bounded from below by a
constant times $(N+K)$, independently of the positions $R_k$ of the nuclei.

By combining the methods in \cite{LY} and our Theorem~\ref{main}, one
can prove the following

\begin{theorem}[{\bf Stability of relativistic matter}]\label{stability}
     There is an $\alpha_{\rm c}>0$ such that for all $N$, $K$,
     $q\alpha\leq\alpha_{\rm c}$ and $\alpha\max\{Z_1,\ldots,Z_K\} \leq 2/\pi$
     one has
     \begin{equation*}
         H_{N,K} \geq -mN.
     \end{equation*}
     The constant $\alpha_{\rm c}$ can be chosen independently of $m$, $A$ and
     $R_1,\ldots,R_K$.
\end{theorem}

The constant $\alpha_{\rm c}$ in Theorem~\ref{stability} depends on the
optimal constant in the Hardy-LT inequality~(\ref{eq:main}) for $d=3$
and $s=1/2$. A bound on this constant, in turn, can be obtained from
our proof in terms of the constant in the Sobolev-Hardy inequalities
(\ref{eq:sobolev}). We do derive a bound on the relevant constants in
Appendix~\ref{app:sobolevconst}, but these bounds are probably far from
optimal. In particular, the available constants do not yield
realistic values of $\alpha_{\rm c}$ so far.

After completing this work, we discovered a different proof of the
special case of the Hardy-LT inequality needed in the proof of
Theorem~\ref{stability}, namely the case where the potential $V$ is
constant inside the unit ball, and infinite outside. In this special
case a substantially improved constant can be obtained, and this
permits the conclusion that Theorem~\ref{stability} holds for the
physical value of $\alpha$, which equals $\alpha\approx 1/137$. We
refer to \cite{FLS} for details.

We briefly outline the proof of Theorem~\ref{stability}. An
examination of the proof in \cite{LY} shows that there are two places
that do not permit the inclusion of a magnetic vector potential $A$.
These are Theorem~9 (Localization of kinetic energy -- general form)
and Theorem~11 (Lower bound to the short-range energy in a ball). Our
Lemma~\ref{imsmag} in Appendix~\ref{app:loc} is precisely the
extension of Theorem~9 to the magnetic case. This lemma implies that
Theorem~10 in \cite{LY} holds also in the magnetic case, without
change except for replacing $|D|$ by $|D-A|$.

Our Theorem~\ref{main} can be used instead of Theorem~11 in \cite{LY}.
In fact, the left side of (3.20) in \cite{LY} is bounded from below by
$q\|\chi\|_\infty^2$ times the sum of the negative eigenvalues of
$|D-A|- \frac 2\pi |x|^{-1} - C R^{-1} \theta_R(x)$, where $\theta_R$
denotes the characteristic function of a ball of radius $R$. By
Theorem~\ref{main}, this latter sum is bounded from below by a
constant times $C^4 R^{-1}$. The resulting bound is of the same form
as the right side of (3.20), except for the constant. It is this
constant that determines the maximally allowed value of the fine
structure constant, $\alpha_{\rm c}$. The rest of the proof remains
unchanged. Note that, in particular, the Daubechies
inequality~\cite{Dau} remains true also in the presence of a magnetic
field.


\subsection{Outline of the paper}
Before giving the proofs of our main results, we pause to outline the
structure of this paper.
\begin{itemize}
\item[$\bullet$]
    In the next Section~\ref{sec:sobolev}, we
give the proof of the Sobolev-Hardy inequalities in
Theorem~\ref{sobolev} and Corollary~\ref{sobolevcor}.
\item[$\bullet$] In the following Section~\ref{sec:rep} we prove what is 
customarily
     called the \lq\lq ground state representation\rq\rq\ in
     Proposition~\ref{substitution}, except that here the \lq\lq ground
     state\rq\rq\ fails to be an $L^2$ function. Such a representation
     for fractional differential operators does not seem to have appeared
     in the literature before.
   \item[$\bullet$] In Section~\ref{sec:lth} we give the proof of our
     main Theorem~\ref{main} about Hardy-LT inequalities. We first
     consider the non-magnetic case $A=0$.  One of the key ingredients
     is the ground state representation obtained in
     Section~\ref{sec:rep}, which allows us to prove a certain
     contraction property of the heat kernel in some weighted
     $L^1$-spaces. Nash's argument \cite[Sect.~8.15--18]{LLo} then
     allows us to translate the Sobolev-Hardy inequalities in
     Corollary~\ref{sobolevcor} into pointwise bounds on the heat kernel
     in an appropriate weighted $L^p$-space. These bounds lead to the
     Hardy-LT inequalities via the trace formula in
     Proposition~\ref{trace} in the spirit of~\cite{L1}.
\item[$\bullet$] Finally, in Section~\ref{sec:diamag} we derive diamagnetic
     inequalities which will allow us to extend the proof of
     Theorem~\ref{main} to the magnetic case.
\end{itemize}

\section{Sobolev-Hardy Inequalities}\label{sec:sobolev}

Our goal in this section is to prove Theorem~\ref{sobolev} and
Corollary~\ref{sobolevcor}. We start with a short outline of the
structure of the proof.

Our proof is based on the fact that we can control the singularity of
$H_s\psi$ near the origin if we know the singularity of $\psi$ at that
point (cf. Lemma~\ref{singularity}). Theorem~\ref{sobolev} follows
by observing that the $L^q$-norm of a symmetric decreasing
function can be bounded above by integrating the function against
$|x|^{d(1/q-1)}$, see Lemma~\ref{isoperimetric}. Moreover, it is
enough to restrict one's attention to symmetric decreasing functions.
Corollary~\ref{sobolevcor} follows from Theorem~\ref{sobolev} by an
IMS-type localization argument, see Lemma~\ref{ims}.

We present some auxiliary results in the following
Subsection~\ref{subsec:aux}. The next two Subsections~\ref{subsec:sob}
and~\ref{subsec:cor} contain the proofs of Theorem~\ref{sobolev} and
Corollary~\ref{sobolevcor}, respectively.

\subsection{Auxiliary material}\label{subsec:aux}

We start with the following integral representation of the operator
$(-\Delta)^{s}$.

\begin{lemma}\label{freeform}
Let $d\geq 1$ and $0<s<1$. Then for all $u\in H^s(\R^d)$
\begin{equation}\label{eq:freeform}
\int_{\R^d} |\xi|^{2s} |\hat u(\xi)|^2 \,d\xi =
a_{s,d} \int_{\R^d} \int_{\R^d} \frac{|u(x)-u(y)|^2}{|x-y|^{d+2s}} \,dx\,dy,
\end{equation}
where
\begin{equation}\label{eq:freeformconst}
a_{s,d} := 2^{2s-1} \pi^{-d/2} \frac{\Gamma((d+2s)/2)}{|\Gamma(-s)|}.
\end{equation}
\end{lemma}

Lemma~\ref{freeform} is well known; we sketch the proof for the sake of
completeness.

\begin{proof}
For fixed $y$ we change coordinates $z=x-y$ and apply Plancherel.
Recalling that $(u(\cdot+z))^\wedge (\xi)=e^{i\xi\cdot z} \hat u(\xi)$
we obtain
\begin{equation*}
\iint \frac{|u(x)-u(y)|^2}{|x-y|^{d+2s}} \,dx\,dy
= \int \left(\int |z|^{-d-2s} |e^{i \xi\cdot z} -1|^2 \,dz \right)
|\hat u(\xi)|^2 \,d\xi\,.
\end{equation*}
The integral in brackets is of the form $c_{s,d} |\xi|^{2s}$, with
\begin{equation*}
\begin{split}
c_{s,d} & := \int_0^\infty \int_{\Sph^{d-1}} |e^{i r\omega\cdot\theta} -
1|^2 \,d\theta \, r^{-2s-1}\,dr \\
& = 2 \int_0^\infty \left( |\Sph^{d-1}| - (2\pi)^{d/2}
r^{-(d-2)/2}J_{(d-2)/2}(r) \right) r^{-2s-1}\,dr\,.
\end{split}
\end{equation*}
Here, $J_{(d-2)/2}$ is the Bessel function of the first kind of order
$(d-2)/2$ \cite{AbSt}. Recall that $|\Sph^{d-1}| = 2\pi^{d/2}/\Gamma(d/2)$. The
formula \eqref{eq:freeformconst} for $c_{s,d}= a_{s,d}^{-1}$ follows now
from
\begin{equation*}
\int_0^\infty r^{-z}\left(J_{(d-2)/2}(r) - 2^{-(d-2)/2}\Gamma(d/2)^{-1}
r^{(d-2)/2}\right)\,dr
= 2^{-z} \frac{\Gamma((d-2z)/4)}{\Gamma((d+2z)/4)}
\end{equation*}
for $d/2<\re z <(d+4)/2$, see \cite[(2.20)]{Ya}.
\end{proof}

Let us recall that $|x|^{-\alpha}$ is a tempered distribution for
$0<\alpha<d$ with Fourier transform
\begin{equation}\label{ft1x}
     b_\alpha \left(|\cdot|^{-\alpha}\right)^\wedge (\xi)
     = b_{d-\alpha} |\xi|^{-d+\alpha},
     \qquad b_\alpha := 2^{\alpha/2} \Gamma(\alpha/2)
\end{equation}
(see, e.g., \cite[Thm.~5.9]{LLo}, where another convention for the
Fourier transform is used, however). We assume now that $s<d/2$. Then
$(-\Delta)^s |x|^{-\alpha}$ is an $L^1_{\rm loc}$-function for
$0<\alpha<d-2s$ and
\begin{equation}\label{eq:powerfunction}
     \left((-\Delta)^s - \mathcal C_{s,d}|x|^{-2s}\right) |x|^{-\alpha}
     = \Phi_{s,d}(\alpha) |x|^{-\alpha-2s}\,,
\end{equation}
where $\mathcal C_{s,d}$ is defined in (\ref{eq:csd}) and
\begin{equation}\label{eq:phi}
\begin{split}
\Phi_{s,d}(\alpha) & := \frac{b_{\alpha+2s}
b_{d-\alpha}}{b_{d-\alpha-2s} b_\alpha} - \mathcal C_{s,d} \\
& = 2^{2s}\left(
\frac{\Gamma((\alpha+2s)/2)\,\Gamma((d-\alpha)/2)}{\Gamma((d-\alpha-2s)/2)\,\Gamma(\alpha/2)}
-\frac{\Gamma^2((d+2s)/4)}{\Gamma^2((d-2s)/4)} \right).
\end{split}
\end{equation}
Later on we will need the following information about the $\alpha$-dependence 
of
$\Phi_{s,d}$.

\begin{lemma}\label{gammamono}
The function $\Phi_{s,d}$ is negative and strictly increasing in
$(0,(d-2s)/2)$ with $\Phi_{s,d}((d-2s)/2)=0$.
\end{lemma}

\begin{proof}
First one checks that
\begin{equation}\label{eq:gammamono1}
\lim_{\alpha\to 0} \Phi_{s,d}(\alpha) = - \mathcal C_{s,d} < 0,
\qquad
\Phi_{s,d}((d-2s)/2)=0.
\end{equation}
Now we abbreviate $\beta:=\alpha/2$, $r:=d/2$ and write
\begin{equation*}
f(\beta) := \Gamma(\beta)/\Gamma(r-\beta),
\qquad
g(\beta) := f(\beta+s)/f(\beta),
\end{equation*}
so that $\Phi_{s,d}(\alpha)=2^{2s} g(\beta)- \mathcal C_{s,d}$. In view
of \eqref{eq:gammamono1} it suffices to verify that $g(\beta)$ is
strictly increasing with respect to $\beta\in (0,(r-s)/2)$. One finds that
\begin{equation}\label{eq:gammamono2}
\frac{f'(\beta)}{f(\beta)}= \frac{\Gamma'(\beta)}{\Gamma(\beta)} +
\frac{\Gamma'(r-\beta)}{\Gamma(r-\beta)}
=\psi(\beta)+\psi(r-\beta)
\end{equation}
with $\psi=\Gamma'/\Gamma$ the Digamma function. Hence
\begin{equation*}
g'(\beta) =
g(\beta) \left(\frac{f'(\beta+s)}{f(\beta+s)} -
\frac{f'(\beta)}{f(\beta)}\right)
=g(\beta) \int_\beta^{\beta+s} h(t)\,dt
\end{equation*}
where, in view of \eqref{eq:gammamono2}, $h(t):= \psi'(t)-\psi'(r-t)$.
Since $\psi'$ is strictly decreasing (see \cite[(6.4.1)]{AbSt}), one has
$h(t)>0$ for $t\in (0,r/2)$. This proves that $g'(\beta)>0$ for all
$\beta\in (0,(r-2s)/2)$. In the case $\beta\in
((r-2s)/2,(r-s)/2)$ one uses in addition the symmetry $h(t)=-h(r-t)$.
\end{proof}


\subsection{Proof of Theorem \ref{sobolev}}\label{subsec:sob}

Our proof of Theorem \ref{sobolev} is close in spirit to \cite{BL} where
a remainder term in the Sobolev inequality on
bounded domains was found. We first exhibit functions $\psi\in\dom h_s$ in
the
form domain which do not lie in the operator domain but for which the
singularity of the distribution (indeed, function) $H_s\psi$ at $x=0$
can be calculated explicitly.

\begin{lemma}\label{singularity}
    Let $0\leq\chi\leq 1$ be a smooth function on $\R_+$ of compact support,
with  $\chi(r)=1$ for $r\leq 1$. Define
\begin{equation}
\psi_\lambda(x) := \chi(|x|/\lambda) |x|^{-\alpha}
\end{equation}
for $0<\alpha<(d-2s)/2$ and $\lambda>0$. Then  $\psi_\lambda \in\dom
h_s$  for  $0<s<1$ and, for every $\epsilon>0$, there exists a
$\lambda_\eps = \lambda_\eps(\alpha,d,s,\chi)$ such that for any
$\lambda\geq \lambda_\eps$,
     \begin{equation}\label{eq:singularity}
       \left(\left((-\Delta)^s -\mathcal C_{s,d} |x|^{-2s}\right)
       \psi_\lambda\right)(x)
       \leq   -\big( |\Phi_{s,d}(\alpha)| -\epsilon\big) |x|^{-\alpha-2s}
       \qquad {\rm for\ all\ } x\in\mathcal B,
     \end{equation}
     in the sense of distributions. Here, $\Phi_{s,d}(\alpha)$ is given in
(\ref{eq:phi}), and $\mathcal B$ denotes the unit ball in $\R^d$.
\end{lemma}

\begin{proof}
     It is not difficult to show that $\psi_\lambda\in H^s(\R^d)$, which implies
     that $\psi_\lambda\in\dom h_s$. (Consult the proof of
     Proposition~\ref{substitution} for details.)  Let $0\leq\phi\in
     C_0^\infty(\mathcal B)$. According to \eqref{eq:powerfunction} one
     has
     \begin{equation*}
       \left(\psi_\lambda,\left((-\Delta)^s -\mathcal C_{s,d}
       |x|^{-2s}\right)\phi\right)
       =
       \Phi_{s,d}(\alpha) \left(|x|^{-\alpha-2s},\phi\right)
       - \left(\tilde\psi_\lambda,(-\Delta)^s\phi\right)\,,
     \end{equation*}
     where $\tilde\psi_\lambda(x):=(1-\chi(|x|/\lambda))|x|^{-\alpha}$. It 
follows from
     Lemma \ref{freeform} (with the aid of polarization) that
     \begin{equation*}
       \left(\tilde\psi_\lambda,(-\Delta)^s\phi\right)
       = -2a_{s,d}\iint 
\frac{\tilde\psi_\lambda(y)\phi(x)}{|x-y|^{d+2s}}\,dx\,dy
       \geq - \rho(\lambda) \int_{\R^d} \frac{
         \phi(x)}{|x|^{\alpha+2s}}\, dx\,,
     \end{equation*}
with
\begin{equation}\label{eq:deltal}
\rho(\lambda)  = \sup_{|x|\leq 1} 2 a_{s,d} |x|^{\alpha+2s} \int
\frac{\tilde\psi_\lambda(y)}{|x-y|^{d+2s}}dy
= \sup_{|x|\leq 1/\lambda} 2 a_{s,d} |x|^{\alpha+2s} \int
\frac{1-\chi(y)}{|y|^\alpha |x-y|^{d+2s}}dy
\,. \end{equation}
Note that $\rho(\lambda)$ is finite for $\lambda\geq 1$, and monotone
decreasing to 0 as  $\lambda\to \infty$. Hence, for a given
$\epsilon>0$ we can choose $\lambda_\epsilon$ such that
$\rho(\lambda_\eps)=\eps$.
Since $\Phi_{s,d}(\alpha)$ is negative by Lemma \ref{gammamono} we
     have established \eqref{eq:singularity}.
\end{proof}

\begin{lemma}\label{isoperimetric}
           Let $1\leq q<\infty$ and $u\in L^q(\R^d)$ a symmetric
decreasing function. Then
           \begin{equation}\label{eq:isoperimetric}
                   \|u\|_q \leq q^{-1}|\mathcal B|^{-1/q'} \int_{\R^d} u(x)
|x|^{-d/q'}\,dx
           \end{equation}
           where $|\mathcal B|$ is the volume of the unit ball $\mathcal B$ in 
$\R^d$,
and $1/q+1/q'=1$.
\end{lemma}

\begin{proof}
           First note that \eqref{eq:isoperimetric} is true (with equality)
if $u$
is the characteristic function of a centered ball. For general $u$ we use
the layer cake representation \cite[Thm.~1.13]{LLo}, $u(x)=\int_0^\infty
\chi_t(x)\,dt$, where $\chi_t$ is the characteristic function of a
centered ball of a certain $t$-dependent radius. Then, by Minkowski's
inequality
\cite[Thm.~2.4]{LLo},
           \begin{align*}
                   \|u\|_q
                   & \leq \int_0^\infty \|\chi_t\|_q \,dt
                   = q^{-1}|\mathcal B|^{-1/q'} \int_0^\infty \int \chi_t(x)
|x|^{-d/q'}\,dx \,dt \\
                   & = q^{-1}|\mathcal B|^{-1/q'} \int u(x) |x|^{-d/q'}\,dx,
           \end{align*}
           proving \eqref{eq:isoperimetric}.
\end{proof}

Now we give the

\begin{proof}[Proof of Theorem \ref{sobolev}]
   We remark first that we may assume $\Omega$ to be a ball and $u$ to
   be a spherically symmetric decreasing function. Indeed, passing to
   the symmetric decreasing rearrangement of $u$ leaves the left side
   of \eqref{eq:sobolev} invariant while it decreases the right side.
   The kinetic energy term on the right side is decreased by virtue of
   Riesz's rearrangement inequality (compare with \cite[Thm.~3.7,
   Lemma~7.17]{LLo}) and $\int |u|^2 |x|^{-2s}\,dx$ increases
   \cite[Thm.~3.4]{LLo}. Moreover, by scaling we may assume that
   $\Omega=\mathcal B$, the unit ball.

Since $\mathcal B$ is bounded, H\"older's inequality implies that it suffices 
to
prove \eqref{eq:sobolev} for $d/(d-2s)<q<2d/(d-2s)$. For such $q$ let
$\alpha:=d/q'-2s$ and note that $0<\alpha<(d-2s)/2$. It follows from Lemmas
\ref{isoperimetric} and \ref{singularity} that for symmetric
decreasing functions $u$ on $\mathcal B$ \begin{equation*}
           \|u\|_q
           \leq q^{-1} |\mathcal B|^{-1/q'} \int_{\mathcal B} u(x)
|x|^{-d/q'} \,dx
     \leq 2 q^{-1} |\mathcal B|^{-1/q'}|\Phi_{s,d}(\alpha)|^{-1} |(u, H_s\psi)|.
\end{equation*}
Here $\psi=\psi_{\lambda_\eps}$ is chosen as in Lemma
\ref{singularity}, with $\eps = |\Phi_{s,d}(\alpha)|/2$. An application of 
Schwarz's
inequality, $|(u, H_s\psi)|^2\leq h_s[u]
h_s[\psi]$, concludes the proof of \eqref{eq:sobolev}.
\end{proof}

In Appendix~\ref{app:sobolevconst}, we shall give an upper bound on the
constant appearing in the Sobolev inequality in the special case of
$d=3$ and $s=1/2$, which is the case of interest in the application in
Section~\ref{sec:stability}.


\subsection{Proof of Corollary \ref{sobolevcor}}\label{subsec:cor}

We will deduce Corollary \ref{sobolevcor} from Theorem~\ref{sobolev}
by a localization argument. For comparison, we recall first the IMS
formula in the local case $s=1$. If $\chi_0,\ldots,\chi_n$ are
Lipschitz continuous functions on $\R^d$ satisfying $\sum_{j=0}^n
\chi_j^2 \equiv 1$, then
\begin{equation}\label{eq:usualims}
     \begin{split}
       & \int_{\R^d}
       \left(|\nabla u|^2-\frac{(d-2)^2}{4}\frac{|u|^2}{|x|^2}
       \right)\,dx \\
       & \qquad = \sum_{j=0}^n \int_{\R^d}
       \left(|\nabla (\chi_j u)|^2-\frac{(d-2)^2}{4}\frac{|\chi_j
       u|^2}{|x|^2} \right)\,dx
       - \int_{\R^d} \sum_{j=0}^n |\nabla \chi_j|^2 |u|^2 \,dx.
     \end{split}
\end{equation}
The following analogous formula for the non-local case, suggested by
Michael Loss, is given in \cite{LY}. The proof is an immediate
consequence of Lemma \ref{freeform}. For a generalization to the
magnetic case, see Lemma~\ref{imsmag} below.

\begin{lemma}\label{ims}
     Let $0<s<\min\{1,d/2\}$ and let
     $\chi_0,\ldots,\chi_n$ be Lipschitz continuous functions on $\R^d$
     satisfying $\sum_{j=0}^n \chi_j^2 \equiv 1$. Then
     \begin{equation}\label{eq:ims}
       h_s[u] = \sum_{j=0}^n h_s[\chi_j u] - (u,L u),
       \qquad u\in C_0^\infty(\R^d),
     \end{equation}
     where $L$ is the bounded operator with integral kernel
     \begin{equation*}
       L(x,y) := a_{s,d} |x-y|^{-d-2s} \sum_{j=0}^n (\chi_j(x)-\chi_j(y))^2.
     \end{equation*}
\end{lemma}

Let us recall the following (non-critical) Sobolev embedding theorem, which is 
easy
to prove. (Cf., e.g., the proof of~\cite[Thm.~8.5]{LLo}.) If $s<d/2$ and
$2\leq q < 2^*=2d/(d-2s)$ then $H^s(\R^d)\subset L^q(\R^d)$ and
\begin{equation}\label{eq:sobolevstandard}
\|u\|_q^2 \leq S_{q,d,s} \left( \|(-\Delta)^{s/2} u\|^2 + \|u\|^2 \right),
\qquad u\in H^s(\R^d).
\end{equation}
In combination with the localization Lemma~\ref{ims} this allows us to
give the

\begin{proof}[Proof of Corollary \ref{sobolevcor}]
Let $\chi_0,\chi_1$ be smooth functions on $\R^d$ with
$\chi_0^2+\chi_1^2\equiv 1$ such that $\chi_0(x)=0$ if $|x|\geq 1$ and
$\chi_1(x)=0$ if $|x|\leq 1/2$. Let $2\leq q < 2d/(d-2s)$. Then, by
Theorem~\ref{sobolev},
\begin{equation*}
\|\chi_0 u\|_q^2 \leq C_{q,d,s} |\mathcal
B|^{2\left(\frac1q-\frac1{2^*}\right)} h_s[\chi_0 u],
\end{equation*}
and by \eqref{eq:sobolevstandard}
\begin{equation*}
\begin{split}
\|\chi_1 u\|_q^2 & \leq S_{q,d,s} \left( \|(-\Delta)^{s/2}(\chi_1 u)\|^2
+ \|\chi_1 u\|^2 \right)\\
& \leq S_{q,d,s} \left( h_s[\chi_1 u] + (2^{2s}\mathcal C_{s,d} +1)
\|\chi_1 u\|^2 \right).
\end{split}
\end{equation*}
Hence Corollary \ref{sobolevcor} follows from Lemma \ref{ims} noting
that $L$ is a bounded operator.
\end{proof}



\section{Ground State Representation}\label{sec:rep}

Eq.~(\ref{eq:powerfunction}) and Lemma~\ref{gammamono} suggest that
the function $|x|^{-(d-2s)/2}$ is a `generalized ground state' of the
operator $H_s$. Our next goal is to establish a ground state
representation. Let us recall the analogous formula in the `local'
case $s=1$. If $d\geq 3$ and $v(x)= |x|^{(d-2)/2}u(x)$ then
\begin{equation}\label{eq:usualsubstitution}
     \int_{\R^d}
     \left(|\nabla u|^2-\frac{(d-2)^2}{4}\frac{|u|^2}{|x|^2} \right)\,dx
     = \int_{\R^d} |\nabla v|^2 \, \frac{dx}{|x|^{d-2}}.
\end{equation}
The corresponding formula in the non-local case $0<s<1$ is
more complicated but close in spirit. It was derived some years ago
by Michael Loss (unpublished notes) for the relativistic case $s=1/2$ and
$d=3$.

\begin{proposition}[{\bf Ground State Representation}]\label{substitution}
Let $0<s<\min\{1,d/2\}$. If $u\in
C_0^\infty(\R^d\setminus\{0\})$ and $v(x)= |x|^{(d-2s)/2}u(x)$, then
\begin{equation}\label{eq:substitution}
h_s[u] =
a_{s,d} \int_{\R^d} \int_{\R^d} \frac{|v(x)-v(y)|^2}{|x-y|^{d+2s}}
\,\frac{dx}{|x|^{(d-2s)/2}}\frac{dy}{|y|^{(d-2s)/2}} \ ,
\end{equation}
with $a_{s,d} $ given in (\ref{eq:freeformconst}).
\end{proposition}

\begin{proof}
Let $0<\alpha<(d-2s)/2$. We shall prove that if $u\in
C_0^\infty(\R^d\setminus\{0\})$ and $v_\alpha(x):=|x|^\alpha u(x)$ then
\begin{equation}\label{eq:substitutionalpha}
           \begin{split}
           & \int_{\R^d} |\xi|^{2s} |\hat u(\xi)|^2 \,d\xi
           - \left(\mathcal C_{s,d}+\Phi_{s,d}(\alpha)\right) \int_{\R^d}
|x|^{-2s} |u(x)|^2 \,dx \\
           & \qquad = a_{s,d} \int_{\R^d} \int_{\R^d}
\frac{|v_\alpha(x)-v_\alpha(y)|^2}{|x-y|^{d+2s}}
                   \,\frac{dx}{|x|^\alpha}\frac{dy}{|y|^\alpha} \ .
           \end{split}
\end{equation}
The proposition follows by letting $\alpha\to (d-2s)/2$. Indeed, the
constant in front of the second integral on the left side then converges to
$\mathcal C_{s,d}$, according to Lemma~\ref{gammamono}. By splitting the
integral into four regions according to the support of $u$, it is easy to
see that the right side is continuous in $\alpha$ and converges to the
right side of \eqref{eq:substitution}.

For the proof of \eqref{eq:substitutionalpha} we can suppose that the
support of $u$ is in the unit ball.
We shall first prove the equality for mollified versions of
$|x|^{-\alpha}$, namely functions $\omega_n(x) = |x|^{-\alpha} \chi(x/n)$,
where $\chi\in C_0^\infty(\R^d)$ with $\chi(x)=1$ for $|x|\leq 1$.

Let us first show that  $\omega_n\in H^s(\R^d)$.  It is clearly in
$L^2(\R^d)$, hence it suffices to establish that
$(-\Delta)^{s/2}\omega_n\in
     L^2(\R^d)$. According to \cite[Thm.~5.9]{LLo} the Fourier transform
     of $\omega_n$ is given by the convolution of $\widehat \chi$ and
     $|\xi|^{\alpha-d}$. Since $\chi$ is assumed to be smooth, $\widehat
     \chi$ decays faster than any power of $|\xi|$. It is then easy to
     see that $\widehat \omega_n$ decays like $|\xi|^{\alpha-d}$, and hence
     $|\xi|^{s}\widehat \psi \in L^2(\R^d)$.

By polarization in Lemma~\ref{freeform}, we get for any $f$ and $g$ in 
$H^s(\R^d)$,
\begin{equation}\label{polarization}
\int_{\R^d} |\xi|^{2s} {\overline{\widehat {f}(\xi)}} \widehat g(\xi)
\,d\xi =
a_{s,d} \iint_{\R^d\times \R^d}  (\overline{f(x)}-\overline{f(y)}
)(g(x)-g(y))\,
\frac{dx\,dy}{|x-y|^{d+2s}} \,.
\end{equation}
We apply this formula to $g(x) = \omega_n(x)$ and $f(x) = |u(x)|^2 /
\omega_n(x)= |u(x)|^2 |x|^\alpha$. In this case, the right side of
(\ref{polarization}) is given by
\begin{equation}\label{pol1}
    a_{s,d} \iint_{\R^d\times\R^d}  \left(
|u(x)-u(y)|^2 - \left|\frac{u(x)}{\omega_n(x)} -\frac{
u(y)}{\omega_n(y)}\right|^2 \omega_n(x)\omega_n(y)\right)
\frac {dx\,dy}{|x-y|^{d+2s}} \,.
\end{equation}
Note that $u(x)/\omega_n(x) = u(x)|x|^\alpha = v_\alpha(x)$ is independent of 
$n$, and
is a $C_0^\infty$ function since the origin is not in the support of
$u$ by assumption.
By dominated convergence, (\ref{pol1}) converges to
\begin{equation}\label{pol2}
    a_{s,d} \iint_{\R^d\times \R^d}  \left(
|u(x)-u(y)|^2 - \left|v_\alpha(x) - v_\alpha(y)\right|^2
|x|^{-\alpha}|y|^{-\alpha}\right)
\frac {dx\,dy}{|x-y|^{d+2s}}
\end{equation}
as $n\to\infty$.  The left side of (\ref{polarization}) can be written as
(compare with (\ref{ft1x}))
\begin{equation}\label{polarizationl}
(2\pi)^{d/2} \frac{b_{d-\alpha}}{b_\alpha} \iint_{\R^d\times \R^d}
|\xi|^{2s} {\overline{\widehat {f}(\xi)}} n^d \widehat
\chi(n(\xi-\xi')) |\xi'|^{\alpha-d}\, d\xi\, d\xi' \,.
\end{equation}
Since $\widehat f$ decays faster than polynomially, $|\cdot|^{2s}
\widehat f
\in L^p(\R^d)$ for any $1\leq p < \infty$. Hence its convolution with
the approximate $\delta$-function $(2\pi)^{d/2} n^d \widehat
\chi(n\,\cdot\,)$ converges to $|\cdot|^{2s} \widehat f $ strongly
in any $L^p$, for $1\leq p<\infty$ \cite[Thm.~2.16]{LLo}. Therefore,
    (\ref{polarizationl}) converges to
\begin{equation}\label{eq:polarization2}
\frac{b_{d-\alpha}}{b_\alpha}\int_{\R^d}  |\xi|^{2s+\alpha-d}
{\overline{\widehat {f}(\xi)}}\,d\xi =
\frac{b_{\alpha+2s} b_{d-\alpha}}{b_{d-\alpha-2s}b_{\alpha}}
\int_{\R^d} |u(x)|^2 |x|^{-2s}\,dx \,.
\end{equation}
Here we used  (\ref{ft1x}) again. The equality of \eqref{pol2} and
\eqref{eq:polarization2} proves \eqref{eq:substitutionalpha}.
\end{proof}

\begin{remark}\label{rem:hardy}
     From the representation (\ref{eq:substitution}) we immediately
     recover Hardy's inequalities (\ref{eq:hardy}) in the case $s<
     \min\{1,d/2\}$. Moreover, we see that equality can not be attained
     (for a non-zero function). From (\ref{eq:substitution}), it is also
     easy to see that the constant $\mathcal C_{s,d}$ is sharp. For this,
     consider a sequence of functions $u_n$, supported in $\mathcal B$,
     approximating $|x|^{-(d-2s)/2}$ close to the origin in a suitable sense.
     The right side of (\ref{eq:substitution}) remains finite in the
     limit $n\to \infty$, whereas $\int |u_n(x)|^2 |x|^{-2s} dx$ diverges.
\end{remark}


\section{Proof of the Hardy-Lieb-Thirring Inequalities}\label{sec:lth}

This section contains the proof of our main result in
Theorem~\ref{main}. We consider here only the non-magnetic case $A=0$,
the extension to non-zero $A$ will be straightforward given the necessary
diamagnetic inequalities which we derive in the next Section~\ref{sec:diamag}.
We explain the necessary modifications in the proof of
Theorem~\ref{main} in Subsection~\ref{sec:dia}.

The ground state representation~(\ref{eq:substitution}) suggests that
it is more natural to regard $h_s[u]$ as a function of $v$ given by
$v(x)= |x|^{(d-2s)/2}u(x)$. In terms of this function $v$, the
Sobolev-Hardy inequality in (\ref{eq:sobolevcor}) can be formulated in
the weighted space with measure $|x|^{-\beta} dx$, where $\beta=
q(d-2s)/2$. Namely,
\begin{equation}\label{eq:Bsob}
    \|v\|^2_{L^q(\R^d,|x|^{-\beta}dx)} \leq C_{q,d,s}'
( v, B_\beta v)_{L^2(\R^d,|x|^{-\beta}dx)}\,,
\end{equation}
where $B_\beta$ is the operator on
$L^2(\R^d,|x|^{-\beta}dx)$ defined by the quadratic form
\begin{equation}\label{defB}
( v, B_\beta v)_{L^2(\R^d,|x|^{-\beta}dx)}=
h_s[ |x|^{-(d-2s)/2} v] + \||x|^{-(d-2s)2}v\|^2_{L^2(\R^d,dx)}\,.
\end{equation}
We suppress the dependence on $s$ in $B_\beta$ for simplicity. Note
that the right side of (\ref{defB}) is independent of $\beta$. The
dependence of $B_\beta$ on $\beta$ comes from the measure
$|x|^{-\beta}dx$ of the underlying $L^2$ space, which is determined by
the value of $q$ in the Sobolev inequality (\ref{eq:Bsob}) as
$\beta=q(d-2s)/2$. We emphasize again that the choice of the weight 
$|x|^{-(d-2s)/2}$ on the right
side of (\ref{defB}) is determined by the ground state representation
(\ref{eq:substitution}).

The proof of Theorem~\ref{main} proceeds in the following steps. 
\begin{itemize}
\item[$\bullet$] From the ground state representation
(\ref{eq:substitution}), we will deduce that $B_\beta$ satisfies the 
Beurling-Deny
criteria, which implies that $e^{-t B_\beta}$ is a contraction on
$L^1(\R^d,|x|^{-\beta} dx)$ (Subsection~\ref{sec:contraction}).
\item[$\bullet$]
Together with the Sobolev-Hardy
inequality (\ref{eq:Bsob}) this yields a bound on the kernel of $e^{-t
     B_\beta}$ via Nash's method (Subsection~\ref{sec:nash}).
\item[$\bullet$] This bound on the heat kernel can then be translated into a
LT bound in the spirit of \cite{L1} (Subsection~\ref{sec:proofmain}).
\end{itemize}

\subsection{Contraction property of $B_\beta$}\label{sec:contraction}

Let $\mathfrak H_\beta :=
L^2(\R^d,|x|^{-\beta}dx)$. We assume that $d-2s < \beta < d$, which
corresponds to $2< q < 2^*$ in (\ref{eq:Bsob}).  The quadratic form
\begin{equation*}
b_\beta[ v ] := h_s[ |x|^{-(d-2s)/2} v ] + \||x|^{-(d-2s)/2} v\|^2\,,
\end{equation*}
considered in the Hilbert space $\mathfrak H_\beta$,
is non-negative and closable on $C_0^\infty(\R^d\setminus\{0\})$, and hence
   generates a self-adjoint operator $B_\beta$ in $\mathfrak H_\beta$.

We shall deduce some positivity properties of the
operator $\exp(-t B_\beta)$.
By Proposition \ref{substitution} the
quadratic form $b_\beta$ satisfies
\begin{enumerate}
\item
if $v,w\in\dom b_\beta$ are real-valued then
$b_\beta[v+iw]=b_\beta[v]+b_\beta[w]$,
\item
if $v\in\dom b_\beta$ is real-valued then $|v|\in\dom b_\beta$ and
$b_\beta[|v|]\leq b_\beta[v]$,
\item
if $v\in\dom b_\beta$ is non-negative then $\min(v,1)\in\dom b_\beta$
and $b_\beta[\min(v,1)]\leq b_\beta[v]$.
\end{enumerate}
By a theorem of Beurling-Deny (see \cite[Section 1.3]{D} or \cite[Section
XIII.12]{ReSi2}) this implies that $\exp(-tB_\beta)$ is
positivity-preserving and a contraction in $L^1(\R^d,
|x|^{-\beta}dx)$. That is, it maps non-negative functions into
non-negative functions, and it decreases $L^1$-norms.

\subsection{Heat kernel estimate}\label{sec:nash}

{}From the contraction property derived above, we will deduce a
pointwise bound on the heat kernel, i.e., on the kernel of the
integral operator $\exp(-t B_\beta)$. We emphasize that this kernel is defined 
by
\begin{equation*}
           \big(\exp(-tB_\beta)v\big)(x) = \int_{\R^d} \exp(-tB_\beta)(x,y)v(y)
\frac{dy}{|y|^\beta}\,.
\end{equation*}
We shall use the Sobolev inequality
\eqref{eq:Bsob} for this purpose.

\begin{proposition}\label{heat}
     Let $d-2s <\beta< d$. Then $\exp(-t B_\beta)$ is an integral
operator on $\mathfrak H_\beta$ and its kernel satisfies
    \begin{equation}\label{eq:heat}
0 \leq \exp(-t B_\beta)(x,y) \leq K_{\beta,d,s} t^{-p} \qquad t>0,\, 
x,y\in\R^d\,,
\end{equation}
     where $p := \beta/(\beta -d + 2s)$. The constant can be chosen to be
$K_{\beta,d,s} := (p\, C_{q,d,s}')^p$ where $C_{q,d,s}'$ is the
constant from Corollary \ref{sobolevcor} with
$q:=2\beta/(d-2s)$.
\end{proposition}

\begin{proof}
Let $\theta:=(q-2)/(q-1)\in (0,1)$. Then H\"older's inequality and
Corollary \ref{sobolevcor} yield for any $v\in
C_0^\infty(\R^d\setminus\{0\})$
\begin{equation*}
\|  v\|^2_{L^2(\R^d,|x|^{-\beta}dx)} \leq \|  v 
\|_{L^q(\R^d,|x|^{-\beta}dx)}^{1-\theta} \| 
v\|_{L^1(\R^d,|x|^{-\beta}dx)}^\theta
\leq C_{q,d,s}'^{(1-\theta)/2}
b_\beta[v]^{(1-\theta)/2} \| v \|_{L^1(\R^d,|x|^{-\beta}dx)}^\theta.
\end{equation*}
Equivalently, if $p$ is as in the proposition, then
\begin{equation}\label{eq:weightednash}
\|v\|_{L^2(|x|^{-\beta}dx)}^{1+1/p} \leq C_{q,d,s}'^{1/2} b_\beta[v]^{1/2}
\|v\|_{L^1(|x|^{-\beta}dx)}^{1/p}.
\end{equation}
This is a Nash-type inequality in $\R^d$ with measure $|x|^{-\beta}dx$.
By Nash's argument (see \cite[Theorem 8.16]{LLo} or \cite[Section
2.4]{D}) this implies that $\exp(-t B_\beta)$ is an integral operator with
kernel satisfying \eqref{eq:heat}, with the constant $K_{\beta,d,s}$
given in the proposition. For the sake of completeness we sketch the
proof of this claim in Appendix \ref{app:nash} below.
\end{proof}

\begin{remark}
    In the above argument we used the Nash-type inequality
    \eqref{eq:weightednash} which we had deduced from the Sobolev-type
    inequality \eqref{eq:sobolevcor}. The heat kernel bound
    \eqref{eq:heat} would actually follow directly from the latter
    inequality by \cite[Thm.~2.4.2]{D}. However, we preferred the
    simplicity of the above argument, yielding in addition an explicit
    constant.
\end{remark}


\subsection{Proof of Theorem \ref{main}}\label{sec:proofmain}

\emph{Step 1.} As a first step, we seek an upper bound on the number of 
eigenvalues
below $-\tau$ of the operator $H_{s}-V$, which we denote by
$N(-\tau,H_{s}-V)$.  By the variational principle we may assume that
$V\geq 0$. Then the Birman-Schwinger principle (see, e.g.,
\cite{ReSi2}) implies that for any increasing non-negative function $F$ on
$(0,\infty)$
\begin{equation*}
N(-1,H_{s}-V)\leq F(1)^{-1} \tr F\left(V^{1/2}(H_{s}+I)^{-1}V^{1/2}\right).
\end{equation*}
Let $\mathcal U : L_2(\R^d) \to \mathfrak H_\beta$ be the unitary
operator which maps $u\mapsto
|x|^{\beta/2}u$. Then
\begin{equation}\label{eq:bs}
     V^{1/2}(H_{s}+I)^{-1}V^{1/2}
     = \mathcal U^* W_\beta^{1/2} B_\beta^{-1} W_\beta^{1/2} \mathcal U \,,
\end{equation}
where $W_\beta$ is the multiplication operator on $\mathfrak H_\beta$
which multiplies by the function $W_\beta(x):=|x|^{\beta+2s-d} V(x)$. Therefore
\begin{equation}
\tr F\left(V^{1/2}(H_{s}+I)^{-1}V^{1/2}\right) = \tr_{\mathfrak H_\beta} 
F\left(W_\beta^{1/2}B_\beta^{-1}W_\beta^{1/2}\right)\,.
\end{equation}
We need the following trace estimate.

\begin{proposition}\label{trace}
     Let $f$ be a non-negative convex function on $[0,\infty)$, growing
polynomially at infinity and
vanishing near the origin, and let
     \begin{equation}\label{eq:f}
       F(\lambda):=\int_0^\infty f(\mu) e^{-\mu/\lambda}\mu^{-1}\,d\mu,
       \qquad \lambda >0.
     \end{equation}
     Then for any $d-2s<\beta<d$ and any multiplication operator $W\geq 0$
     \begin{equation}\label{eq:trace}
      \tr_{\mathfrak H_\beta} F\left(W^{1/2}B_\beta^{-1}W^{1/2}\right)
       \leq \int_0^\infty\int_{\R^d} \exp(-t B_\beta)(x,x)
       f(t W(x))\,\frac{dx}{|x|^{\beta}} \frac{dt}t.
     \end{equation}
\end{proposition}

Note that the heat kernel $\exp(-tB_\beta)(x,y)$ is well defined on the
diagonal $x=y$ by the semigroup property. Namely, $\exp(-tB_\beta)(x,x)
= \int |\exp(-t B_\beta/2)(x,y)|^2 |y|^{-\beta} \,dy$. For the proof of 
Proposition~\ref{trace} one follows the proof of the CLR
bound in \cite{L1} (see also \cite{Si2} and \cite{RoSo}). As in the
latter paper Trotter's product formula can be used in place of path
integrals. For details we refer to Appendix~\ref{app:trace}.

We shall now assume that $F$ has the special form \eqref{eq:f} in
order to apply the trace estimate from Proposition \ref{trace}. Given
$d-2s < \beta<d$ and $p= \beta/(\beta -d + 2s)$,
Proposition~\ref{heat} implies that
           \begin{equation*}
           \begin{split}
             & \int_0^\infty\int_{\R_d} \exp(-t B_\beta)(x,x)
             f(t W_\beta(x))\, \frac{dx}{|x|^\beta} \frac{dt}t \\
             & \qquad \leq K_{\beta,d,s} \int_0^\infty \int_{\R^d} t^{-p}  f(t 
W_\beta(x))\,\frac{dx}{|x|^\beta}\frac{dt}t \\
             & \qquad = K_{\beta,d,s}\left(\int_{\R^d} W_\beta(x)^p
               \frac{dx}{|x|^\beta} \right)
             \left(\int_0^\infty t^{-p-1} f(t)\,dt\right) \,.
           \end{split}
           \end{equation*}
Note that $W_\beta(x)^p= V(x)^p|x|^{\beta}$. We conclude that for any $d/2s < p 
< \infty$, \begin{equation}\label{eq:number}
     N(-1,H_{s} - V) \leq K_{p,d,s}' \int_{\R^d} V(x)_+^p \,dx\,,
     \end{equation}
where the constant is given by \begin{equation*}
K_{p,d,s}' =
K_{\beta,d,s}\, \inf_f F(1)^{-1}
\left(\int_0^\infty t^{-p-1} f(t)\,dt\right)\,.
\end{equation*}
Here $\beta=p(d-2s)/(p-1)$, and the infimum runs over all admissible
functions $f$ from Proposition \ref{trace}. In order to obtain an
explicit upper bound one may choose $f(x):=(x-a)_+$ and minimize
over $a>0$.

\emph{Step 2.} Now we use the idea of \cite{LTh} to deduce
\eqref{eq:main} from \eqref{eq:number}. Fix $\gamma>0$ and choose some
$d/(2s)<p<\gamma+d/(2s)$. First we note that by scaling we have, for any
$\tau>0$,
\begin{equation*}
N(-\tau,H_{s} - V) = N(-1, H_{s} - V_\tau)
\end{equation*}
where $V_\tau(x):=\tau^{-1}V(\tau^{-1/2s}x)$. In view of
\eqref{eq:number} this yields
\begin{equation}\label{eq:numberscale}
N(-\tau,H_{s} - V) \leq K_{p,d,s}'\tau^{-p+d/2s} \int_{\R^d} V(x)_+^p \,dx.
\end{equation}
Now, for any fixed $0<\sigma<1$, one has by the variational principle
\begin{equation*}
N(-\tau,H_s-V) \leq N(-(1-\sigma)\tau, H_s-(V-\sigma \tau)_+)\,.
\end{equation*}
Hence, by \eqref{eq:numberscale},
\begin{align*}
     \tr (H_{s}-V)_-^\gamma
     & = \gamma \int_0^\infty N(-\tau,H_{s}-V) \tau^{\gamma-1}\,d\tau \\
     & \leq \gamma K_{p,d,s}' (1-\sigma)^{-p+d/2s}  \int_0^\infty 
\!\!\int_{\R^d} (V(x)-\sigma \tau)_+^p \,dx
     \, \tau^{\gamma-p+d/2s-1}\,d\tau\,.
\end{align*}
We change the order of integration and calculate the $\tau$-integral
first. For fixed  $x\in\R^d$, \begin{equation*}
       \int_0^\infty (V(x)-\sigma \tau)_+^p \tau^{\gamma-p+d/2s-1}\,d\tau
       = \sigma^{-\gamma-d/2s+p} V(x)_+^{\gamma+d/2s}
       B(\gamma+d/2s-p,p+1)\,.
\end{equation*}
Here, $B$ denotes the Beta-function
$B(a,b)=\Gamma(a)\Gamma(b)/\Gamma(a+b)$.  Minimization over $\sigma\in
(0,1)$ and $p \in (d/2s, \gamma+d/2s)$ yields
\begin{equation}
       \tr (H_{s}-V)_-^\gamma
       \leq C_{d,s}(\gamma) \int_{\R^d} V(x)_+^{\gamma+d/2s}\,dx \label{eq:lth}
\end{equation}
with
\begin{align*}
     C_{d,s}(\gamma) := \min_{d/2s< p <\gamma+d/2s} & \biggl\{ \gamma^{\gamma+1} 
K_{p,d,s}' B(\gamma+d/2s-p,p+1)\\
     & \quad \times (\gamma+d/2s-p)^{-\gamma-d/2s+p} (p-d/2s)^{-p+d/2s} 
\biggl\}\,.
\end{align*}
This concludes
the proof of Theorem \ref{main} in the case $A=0$.


\section{Extension to Magnetic Fields}\label{sec:diamag}

In this section we prove certain diamagnetic inequalities which allow
us to extend the proof of Theorem~\ref{main} in the previous section
to the case of non-zero magnetic fields. The main idea is contained in
Proposition~\ref{diamaggen} in Subsection~\ref{sec:diamag1}. The
following Subsection~\ref{sec:diamag2} contains some technical
refinements we will need. In Subsection~\ref{sec:dia} we describe the
necessary modifications in the proof of Theorem~\ref{main} to include
magnetic fields. The final Subsection~\ref{sec:s1} is devoted to the special
case $s=1$.

Throughout this section we shall assume that $A\in
L^2_{\rm loc}(\R^d;\R^d)$ and that $d\geq 2$. Note that in $d=1$ any
magnetic vector potential can be removed by a gauge transformation.

In $\mathfrak H_\beta= L^2(\R^d,|x|^{-\beta}dx)$ consider the
quadratic form
\begin{equation}\label{def:bba}
b_{\beta,A}[ v ] := h_{s,A}[ |x|^{-(d-2s)/2} v ] + \||x|^{-(d-2s)/2}
v\|^2\quad,\quad v\in C_0^\infty(\R^d\setminus \{0\})\,.
\end{equation}
In Appendix~\ref{app:closable}, we show that $b_{\beta,A}$ is
closable and hence defines a self-adjoint operator
$B_{\beta,A}$ in $\mathfrak H_\beta$. Our goal in this section is to
show that $\exp(-tB_{\beta,A})$ is an integral operator, whose kernel
satisfies
\begin{equation}\label{eq:ba}
|\exp(-tB_{\beta,A})(x,y)| \leq \exp(-tB_\beta)(x,y)\,.
\end{equation}

\subsection{An initial inequality}\label{sec:diamag1}

We consider weighted magnetic operators $\omega |D-A|^{2s} \omega$ where
$A\in L^2_{\rm loc}(\R^d)$ and $\omega>0$ with $\omega+\omega^{-1}\in
L^\infty(\R^d)$. This is a self-adjoint operator in $L^2(\R^d)$ with form
domain $\omega^{-1} \dom|D-A|^s$. It satisfies the following diamagnetic 
inequality.

\begin{proposition}[{\bf Weighted diamagnetic inequality}]\label{diamaggen}
           Let $d\geq 2$ and $0<s\leq 1$. Assume that $A\in L^2_{\rm loc}(\R^d)$ 
and
that $\omega>0$ with $\omega+\omega^{-1}\in L^\infty(\R^d)$. Then for all
$u\in L^2(\R^d)$ and all $t\geq 0$ one has
     \begin{equation}\label{eq:diamaggen}
       |\exp(-t \omega |D-A|^{2s} \omega) u|
       \leq \exp(-t \omega(-\Delta)^s\omega) |u|.
     \end{equation}
\end{proposition}

\begin{proof}
First note that the assertion is true in the case $\omega\equiv
1$, i.e.,
\begin{equation}\label{eq:diamagheat}
|\exp(-t|D-A|^{2s})u| \leq \exp(-t(-\Delta)^s)|u|.
\end{equation}
Indeed, for $s=1$ this inequality is proved in \cite{Si1} for all
$A\in L^2_{\rm loc}(\R^d)$. The general case $0<s<1$ follows from the
fact that the function $\lambda \mapsto e^{-\lambda^s}$ is completely
monotone (i.e., its derivatives are alternating in sign) and hence is the 
Laplace
transform of a positive measure by Bernstein's theorem~\cite{Do}. This
reduces the problem to the case $s=1$.

Now assume that $\omega$ is as in the proposition and write $M_A:=\omega
|D-A|^{2s} \omega$. In view of the general relation
           \begin{equation}\label{eq:expres}
                   \exp(-tM_A) = \slim_{n\to\infty}\, (I+n^{-1}tM_A)^{-n}
           \end{equation}
           it suffices to prove the inequality $|(M_A+\tau)^{-1} u| \leq
(M_0+\tau)^{-1}|u|$. But since $(M_A+\tau)^{-1} = \omega^{-1} (|D-A|^{2s}
+ V)^{-1} \omega^{-1}$ with $V:=\tau \omega^{-2}$ it suffices to prove
$|(|D-A|^{2s} + V)^{-1} u| \leq (-\Delta^{s} + V)^{-1}|u|$. In view of
the relation  `inverse'  to \eqref{eq:expres},
           \begin{equation*}
                   (|D-A|^{2s} + V)^{-1} = \int_0^\infty \exp(-t(|D-A|^{2s} +
V))\,dt,
           \end{equation*}
           the assertion follows from \eqref{eq:diamagheat} and Trotter's
product
formula.
\end{proof}

\begin{remark}\label{diamag1}
           We note that the diamagnetic inequality in the form
\eqref{eq:diamagheat}
implies \eqref{eq:diamag1}. This follows by integrating the square of
\eqref{eq:diamagheat} and evaluating the derivative with respect to $t$
at $t=0$.
\end{remark}


\subsection{A refined inequality}\label{sec:diamag2}

In order to prove the desired inequality (\ref{eq:ba}), we have to extend
Proposition \ref{diamaggen} in two directions. First, we want to use
the singular weight $\omega(x)=|x|^\alpha$, $0<\alpha<s$, and second,
we want to replace the operator $|D-A|^{2s}$ by $H_{s,A}+I$, i.e., we
want to subtract the Hardy term.

Recall that $B_{\beta,A}$ was defined by (\ref{def:bba}). The main
result of this section is the following.

\begin{proposition}[{\bf Weighted diamagnetic inequality, second
version}]\label{diamag}
           Let $d\geq 2$ and $d-2s < \beta < d$. Assume that $A\in
L^2_{\rm loc}(\R^d)$.
Then for all $v\in \mathfrak H_\beta$ one has
     \begin{equation}\label{eq:diamag}
       |\exp(-t B_{\beta,A})v| \leq \exp(-t B_\beta)|v|.
     \end{equation}
\end{proposition}

It follows from Proposition~\ref{heat} that $\exp(-t B_\beta)$ is an integral
operator that maps $L^1(\R^d,|x|^{-\beta}dx)$ to
$L^\infty(\R^d,|x|^{-\beta}dx)$. Hence (\ref{eq:diamag}) implies that the
same is true for $\exp(-t B_{\beta,A})$. Moreover, the kernels are
related by the inequality (\ref{eq:ba}).

In the course of the proof we will need the following approximation result.

\begin{lemma}\label{approx}
           Let $T_n$, $T$ be closed, densely defined operators in a Hilbert
space
$\mathfrak H$ with $T_n T_n^* \to TT^*$ in strong resolvent sense. Assume
that there is a set $\mathcal D\subset\bigcap\dom T_n\cap\dom T$, dense
in $\mathfrak H$, such that $T_n\phi\to T\phi$ for all $\phi\in\mathcal
D$. Then $T_n^* T_n \to T^*T$ in strong resolvent sense.
\end{lemma}

\begin{proof}
           For $\gamma>0$ and $\phi,\psi\in\mathcal D$ one has
           \begin{equation*}
                   \gamma (\phi,(T_n^* T_n+\gamma)^{-1}\psi) = (\phi,\psi)
-(T_n\phi, (T_n
T_n^*+\gamma)^{-1} T_n\psi).
           \end{equation*}
           By assumption the right side converges to $(\phi,\psi) - (T\phi, (T
T^*+\gamma)^{-1} T\psi)= \gamma (\phi,(T^* T+\gamma)^{-1}\psi)$, which
proves that $T_n^* T_n \to T^*T$ in weak resolvent sense. However, the
latter is the same as in strong resolvent sense, see \cite[Problem
VIII.20]{ReSi1}.
\end{proof}

\begin{proof}

\emph{Step 1.} Consider again the unitary transformation $\mathcal U :
L^2(\R^d) \to \mathfrak H_\beta$, which maps $u\mapsto
|x|^{\beta/2}u$. Then \begin{equation}\label{def:qbu}
Q_{\beta,A} := \mathcal U^* B_{\beta,A} \mathcal U \end{equation}
is a self-adjoint operator in the unweighted space $L^2(\R^d)$, whose
quadratic form is given by
\begin{equation}\label{def:qba}
     q_{\beta,A}[u] := h_{s,A}[|x|^\alpha u] + \||x|^\alpha u\|^2 \quad,
     \quad \alpha = (\beta+2s-d)/2 \,.
\end{equation}
Eq. (\ref{eq:diamag}) is equivalent to
     \begin{equation}\label{eq:diamag2}
       |\exp(-t Q_{\beta,A})u| \leq \exp(-t Q_{\beta,0})|u|.
     \end{equation}
for $u\in L^2(\R^d)$.

\emph{Step 2.}  We begin by considering the case where the `potential
terms' in the definition of $Q_{\beta,A}$ are absent. More precisely,
we consider the operator $M_{\beta,A}$ in $L^2(\R^d)$ generated by the
quadratic form $\| |D-A|^{s} |x|^\alpha u\|^2$ on
$C_0^\infty(\R^d\setminus\{0\})$. We shall prove that for all $u\in
L^2(\R^d)$ one has
\begin{equation}\label{eq:diamagpure}
       |\exp(-t M_{\beta,A}) u| \leq \exp(-t M_{\beta,0}) |u|.
\end{equation}

Let $\omega_n$ be a family of smooth positive functions which decrease
monotonically to $|x|^\alpha$ and agree with this function outside a
ball a radius $n^{-1}$. Similarly, for fixed $n$ let $\omega_{n,m}$ be
a family of smooth positive and bounded functions which increase
monotonically to $\omega_n$ and agree with this function inside a ball
a radius $m$.

The operators $|D-A|^s \omega_n$ and $|D-A|^s \omega_{n,m}$ are easily
seen to be closable on $C_0^\infty(\R^d\setminus\{0\})$ (as in
Appendix \ref{app:closable}) and we denote their closures by $T_n$ and
$T_{n,m}$, respectively. One finds that $C_0^\infty(\R^d)\subset\dom
T_n^*$ with $T_n^* v = \omega_n |D-A|^s v$ for $v\in
C_0^\infty(\R^d)$, and similarly for $T_{n,m}^*$. By construction of
$\omega_{n,m}$ the operators $T_{n,m}T_{n,m}^*$ are monotonically
increasing as $m\to\infty$ and hence converge in strong resolvent
sense to $T_{n}T_{n}^*$ by \cite[Thm.~S.14]{ReSi1}. Noting that
$T_{n,m}\phi\to T_n\phi$ for any $\phi\in
C_0^\infty(\R^d\setminus\{0\})$ we conclude from Lemma \ref{approx}
that $T_{n,m}^*T_{n,m}\to T_{n}^* T_{n}$ in strong resolvent sense.
One checks that $T_{n,m}^*T_{n,m}$ coincides with the operator
$\omega_{n,m}|D-A|^{2s}\omega_{n,m}$ from Subsection \ref{sec:diamag1}
and satisfies a diamagnetic inequality by Proposition \ref{diamaggen}.
By the strong resolvent convergence the diamagnetic inequality is also
valid for $T_{n}^*T_{n}$. Now we repeat the argument for $n\to\infty$
where we have monotone convergence from above. We apply
\cite[Thm.~1.2.3]{D} noting that $C_0^\infty(\R^d\setminus\{0\})$ is a
form core for all the operators involved, and conclude that
$T_{n}^*T_{n}\to M_{\alpha,A}$ in strong resolvent sense. This proves
the diamagnetic inequality \eqref{eq:diamagpure}.

\emph{Step 3.}  Now we use another approximation argument to include
the Hardy term. We define $R_{\beta,A}$ via the quadratic form
$h_{s,A}[|x|^\alpha u]$ on $C_0^\infty(\R^d\setminus\{0\})$. Moreover,
for $n\in\N$ let $W_n(x):= \mathcal C_{s,d}
\min\{|x|^{-2(s-\alpha)},n\}$. The boundedness of $W_n$,
\eqref{eq:diamagpure} and Trotter's product formula show that the
diamagnetic inequality is valid for $M_{\beta,A}-W_n$. Since
$M_{\beta,A}-W_n \to R_{\beta,A}$ in strong resolvent sense again by
monotone convergence we find the diamagnetic inequality
\eqref{eq:diamagpure} with $M_{\beta,A}$ replaced by $R_{\beta,A}$.

\emph{Step 4.}  Finally, we note that $Q_{\beta,A}= R_{\beta,A} +
|x|^{2\alpha}$ in the sense of quadratic forms. Moreover,
$C_0^\infty(\R^d\setminus\{0\})$ is a core for both quadratic forms
involved. Equation~(\ref{eq:diamag2}) follows now from the diamagnetic
inequality for $R_{\beta,A}$ by Kato's strong Trotter product formula
\cite[Theorem S.21]{ReSi1}.
\end{proof}

\subsection{Extension of Theorem~\ref{main} to magnetic fields}\label{sec:dia}
In the case of non-vanishing magnetic field, the proof of
Theorem~\ref{main} is essentially identical to the one presented in
the previous section. Although
(\ref{eq:trace}) does not necessarily hold with $B_{\beta,A}$ instead
of $B_\beta$, it \emph{does} hold if  $B_\beta$ is replaced by
$B_{\beta,A}$ on the left side only! I.e.,
\begin{equation*}\label{eq:trace2}
\tr_{\mathfrak H_\beta} F\left(W^{1/2}B_{\beta,A}^{-1}W^{1/2}\right)
    \leq \int_0^\infty\int_{\R^d} \exp(-t B_\beta)(x,x)
     f(t W(x))\,\frac{dx}{|x|^{\beta}} \frac{dt}t\,.
\end{equation*}
For the proof, one uses the diamagnetic inequality (\ref{eq:ba}) before 
applying Jensen's
inequality (cf. Appendix~\ref{app:trace}).
This leads to the conclusion that  (\ref{eq:number}) also
holds with magnetic fields, i.e.,
\begin{equation}\label{eq:number2}
     N(-1,H_{s,A} - V) \leq K_{p,d,s}' \int_{\R^d} V(x)_+^p \,dx\,,
     \end{equation}
with the same ($A$-independent) constant $K_{p,d,s}'$.

For the remainder of the proof, we note that for any
$\tau>0$,
\begin{equation*}
N(-\tau,H_{s,A} - V) = N(-1, H_{s,A_\tau} - V_\tau)
\end{equation*}
where $V_\tau(x):=\tau^{-1}V(\tau^{-1/2s}x)$ and
$A_\tau(x):=\tau^{-1/2s}A(\tau^{-1/2s}x)$. The scaling of $A$ does not
have any effect, however, since the constant in (\ref{eq:number2}) is
independent of $A$. Therefore (\ref{eq:lth}) also holds with $H_s$
replaced by $H_{s,A}$, with same constant $C_{d,s}(\gamma)$.

\subsection{The special case $s=1$}\label{sec:s1}

As noted in Remark~\ref{rem:s1}, the proof of Theorem~\ref{main} just
given works also in the case $s=1$ in dimensions $d\geq 3$. We briefly comment
on the necessary modifications.

The local Sobolev-Hardy inequalities for $s=1$ have been proved in
\cite{BVa}. Alternatively, one can obtain them following our proof in
Section~\ref{sec:sobolev}. Using the IMS formula (\ref{eq:usualims})
one can obtain the global Sobolev-Hardy inequalities
(\ref{frank1}). The rest of the proof goes through without change. To
verify the Beurling-Deny criteria, one uses the ground-state
representation (\ref{eq:usualsubstitution}) instead of
Proposition~\ref{substitution}. Note also that the weighted diamagnetic
inequalities in this section include the case $s=1$.


\begin{appendix}

\section{A Constant in the Sobolev Inequality
\eqref{eq:sobolev}}\label{app:sobolevconst}

In this appendix we shall derive an explicit bound on the constants
$C_{q,3,1/2}$  for the Sobolev-Hardy inequalities (\ref{eq:sobolev}) in the 
case
$d=3$ and $s=1/2$, which is of interest for our theorem on
stability of matter. Let $3/2<q<3$ and $\alpha:=2-3/q$. For $\lambda >1$, let 
\begin{equation}\label{eq:defrho}
\rho(\lambda) := \frac{1-\alpha}{\pi \lambda^{1+\alpha}}
\int_1^\infty \frac{dr}{r^{(1+\alpha)/2} (r-\lambda^{-2})}\,.
\end{equation}
  We will show that \begin{equation}\label{eq:c3}
C_{q,3,1/2} \leq \frac{\pi^2}{3 q^2} \left(1-\alpha \right) (3/4\pi)^{4/3}
  \inf_{\lambda>1}\frac{ \lambda^{2(1-\alpha)}}{ \big(|\Phi_{1/2,3}(\alpha)|- 
\rho(\lambda)\big)_+^{2}} \,.
\end{equation}
We remark
that in this special case
\begin{equation*}\label{eq:phiexpl}
|\Phi_{1/2,3}(\alpha)| = \frac 2\pi -(1-3/q)\, \cot \left(\frac{\pi
(1-3/q)}2\right) \,.
\end{equation*}

The estimate (\ref{eq:c3}) is a consequence of the following two
facts. First, we claim that
\begin{equation}\label{eq:psiestimate}
\limsup_{\delta\to 0}\limsup_{\epsilon \to 0}
h_{1/2}[\psi^{\epsilon,\delta}_\lambda ] \leq
    \lambda^{2(1-\alpha)} \frac{\pi^2}{3}\left( 1-\alpha\right)\,,
\end{equation}
where $\psi^{\epsilon,\delta}_\lambda(x) =
\lambda^{-\alpha}\psi^{\epsilon,\delta}(x/\lambda)$ and 
$\psi^{\epsilon,\delta}$ is defined for $\epsilon,\delta>0$ by
\begin{equation}\label{eq:psi}
\psi^{\epsilon,\delta}(x) = \left\{ \begin{array}{ll}
|x|^{-\alpha} & {\rm for\ } |x|\leq 1\,, \\
|x|^{-1} \left( 1  -  \epsilon^\delta (|x|^2-1)^\delta \right) & {\rm for\ }
1\leq |x|^2 \leq 1 + 1/\epsilon\,, \\
0 & {\rm for\ } |x|^2 \geq 1 + 1/\epsilon\,.
\end{array}\right.
\end{equation}
Note that $\psi^{\eps,\delta}$ does not satisfy the smoothness
assumption of Lemma~\ref{singularity}, but it can be approximated by
such functions in $h_{1/2}$-norm.

Secondly, we claim that $\rho^{\eps,\delta}$ in (\ref{eq:deltal})
defined with the function $\psi_\lambda =
\psi^{\eps,\delta}_\lambda$ satisfies \begin{equation}\label{eq:rhoeps}
\lim_{\eps \to 0} \rho^{\eps,\delta}(\lambda) = \rho(\lambda)
\end{equation}
uniformly in $\delta>0$ and $\lambda>1$, with $\rho(\lambda)$ as in
(\ref{eq:defrho}). Eq. (\ref{eq:c3}) then follows from these two
facts, proceeding as in the proof of Theorem~\ref{sobolev}. Instead of
choosing $\lambda$ such that $\rho(\lambda)= |\Phi_{1/2,3}(\alpha)|/2$
we optimize now over the choice of $\lambda$.

For the proof of \eqref{eq:psiestimate} we consider first an arbitrary
radial function $\psi$ and, with a slight abuse of notation, we write
$\psi(x) = \psi(r)$ for $r=|x|$. Using the ground state representation
in Proposition~\ref{substitution}, introducing spherical coordinates, and
integrating over the angles, we have
\begin{equation*}
h_{1/2}[\psi] = 8 \int_0^\infty dr
\int_0^\infty ds \frac{rs}{(r^2-s^2)^2} \left| r \psi(r) -
     s\psi(s)\right|^2\,.
\end{equation*}
By changing variables $r^2\to r$ and $s^2 \to s$, this yields
\begin{equation}\label{eq:radialinte}
h_{1/2}[\psi] = 2 \int_0^\infty dr
\int_0^\infty ds \frac{1}{(r-s)^2} \left| \sqrt{r} \psi(\sqrt{r}) -
     \sqrt{s}\psi(\sqrt{s})\right|^2\,.
\end{equation}

Now assume that $\psi=\psi^{\epsilon,\delta}$ as in \eqref{eq:psi}. By
scaling, if suffices to prove (\ref{eq:psiestimate}) for $\lambda=1$. We
split the integrals
in (\ref{eq:radialinte}) into several parts. First of all, we have
\begin{equation}\label{eq:app14}
    \int_0^1 dr
\int_0^1 ds \frac{1}{(r-s)^2} \left( r^{(1-\alpha)/2} -
     s^{(1-\alpha)/2}\right)^2 = 2  \int_0^1 ds
    \frac 1{s^\alpha} \int_0^1 dt \frac{1}{(1-t)^2} \left( 1 -
     t^{(1-\alpha)/2}\right)^2 \,.
\end{equation}
This identity can be obtained by noting that the integral on the left
is the same as twice the integral over the region $r\leq s$, and then
writing $r= st$ for $0\leq t\leq 1$. Simple computations then lead to
\begin{align}\nonumber
    \int_0^1 dt \frac{1}{(1-t)^2} \left( 1 -
     t^{(1-\alpha)/2}\right)^2 & = \frac{(1-\alpha)^2}4 \int_0^1 dt
\frac{1}{(1-t)^2} \int_t^1 ds\, s^{-(1+\alpha)/2} \int_t^1 du\,
u^{-(1+\alpha)/2} \\ \nonumber & = \frac{(1-\alpha)^2}2  \int_0^1 ds
\int_s^1 du\, (su)^{-(1+\alpha)/2} \int_0^s dt \frac 1{(1-t)^2} \\ & =
(1-\alpha) \int_0^1 ds \frac 1{1-s} s^{(1-\alpha)/2} \left( 1-
     s^{(1-\alpha)/2}\right) \,.
\end{align}
We introduce the function
\begin{equation*}
\eta(\lambda) = \int_0^\infty dt \left( \frac{e^{-t}}t -
     \frac{e^{-\lambda t}}{1-e^{-t}}\right) \,.
\end{equation*}
We note that $\eta(\lambda)= \Gamma'(\lambda)/\Gamma(\lambda)$ is the
Digamma-function. It is then easy to see that
\begin{equation}\label{eq:app17}
    \int_0^1 ds \frac 1{1-s} s^{(1-\alpha)/2} \left( 1-
     s^{(1-\alpha)/2}\right) = \eta(2 -\alpha) - \eta(3/2 - \alpha/2) \,.
\end{equation}
Altogether, we conclude that the contribution of $r\leq 1$ and $s\leq
1$ to the integral in (\ref{eq:radialinte}) is given by
\begin{equation*}
4 \left(\eta(2-\alpha) - \eta(3/2 -\alpha/2) \right) \,.
\end{equation*}

Similarly, we proceed with the other terms. We have
\begin{align}\nonumber
&\lim_{\eps \to 0}  \int_0^1 dr
\int_1^{\infty} ds \frac{1}{(r-s)^2} \left( r^{(1-\alpha)/2} -
     \left[ 1 -  \eps^\delta (s-1)^\delta\right]_+ \right)^2 \\  &=  \int_0^1 dr
\int_1^{\infty} ds \frac{1}{(r-s)^2} \left( r^{(1-\alpha)/2} - 1
     \right)^2 =  \int_0^1 dr
\frac{1}{1-r} \left( r^{(1-\alpha)/2} - 1
     \right)^2\,.
\end{align}
Here we have used dominated convergence, noting that the integrand is
bounded from above by the $L^1$ function $(r-s)^{-2}( 2
(1-r^{(1-\alpha)/2})^2+2\min\{1,(s-1)^{2\delta}\})$ for $\eps\leq 1$.
The contribution of this term to (\ref{eq:radialinte}) (noting that it
appears twice) is thus given
by
\begin{equation*}
4 \left( 2 \eta(3/2-\alpha/2) - \eta(1) - \eta(2-\alpha)\right)\,.
\end{equation*}

We are left with calculating
\begin{align*}\nonumber
&  \int_1^{1+1/\eps} dr
\int_1^{1+1/\eps} ds \frac{1}{(r-s)^2} \left( \eps^\delta (r-1)^\delta -
     \eps^\delta (s-1)^\delta\right)^2 \\  &= \int_0^{1} dr
\int_0^{1} ds \frac{1}{(r-s)^2} \left( r^\delta -
     s^\delta\right)^2 = 2 \left( \eta(1+2\delta) - \eta(1+\delta)\right) \,.
\end{align*}
The last equality follows by proceeding as in
(\ref{eq:app14})--(\ref{eq:app17}). The last term to evaluate is
\begin{align*}\nonumber
&  \int_1^{1+1/\eps} dr
\int_{1+1/\eps}^\infty ds \frac{1}{(r-s)^2} \left( 1- \eps^\delta
     (r-1)^\delta\right)^2 \\  &=  \int_0^{1} dr
    \frac{1}{1-r} \left( 1- r^\delta\right)^2 = 2 \eta(1+\delta) -
    \eta(1) - \eta(1+2\delta)\,.
\end{align*}

We have thus shown that
\begin{equation*}
\lim_{\epsilon \to 0} h_{1/2}[\psi^{\epsilon,\delta}] =
4 \left( \eta(3/2 -\alpha/2) + \eta(1+\delta) - 2 \eta(1) \right) \,.
\end{equation*}
Using concavity of $\eta$, together with $\eta'(1) = \pi^2/6$, yields
the estimate
\begin{equation*}
\lim_{\epsilon \to 0} h_{1/2}[\psi^{\epsilon,\delta}] \leq
     \frac{2 \pi^2}{3}\left( \frac{ 1-\alpha}2 + \delta\right)\,.
\end{equation*}

We proceed similarly for the calculation of
$\rho^{\eps,\delta}$. We have $$
\rho^{\eps,\delta}(\lambda) = \frac{2}{\pi \lambda^{1+\alpha}}
\int_1^{1+1/\eps} \frac{ r^{(1-\alpha)/2}}{(r-\lambda^{-2})^2} \left(
   1- \frac{1-\eps^\delta(r-1)^\delta}{r^{(1-\alpha)/2}}\right)dr\,.
$$
Eq. (\ref{eq:rhoeps}) then follows by dominated convergence and
integration by parts.


\section{Localization Formula in the Magnetic Case}\label{app:loc}

In this appendix we establish the analogue of Proposition \ref{ims} in
the general case $A\not=0$. As explained in
Section~\ref{sec:stability}, this is needed for the proof of
Theorem~\ref{stability}. First recall that for $s=1$ and
$\sum_{j=1}^n\chi_j^2\equiv 1$ one has
\begin{equation*}
     \int_{\R^d} |(D-A)u|^2 \,dx
     = \sum_{j=0}^n \int_{\R^d}
     |(D-A) (\chi_j u)|^2 \,dx
     - \int_{\R^d} \sum_{j=0}^n |\nabla \chi_j|^2 |u|^2 \,dx.
\end{equation*}
In this case the localization error $\sum_{j=0}^n |\nabla \chi_j|^2$ is
local and independent of $A$. The analogue for $s<1$ is

\begin{lemma}\label{imsmag}
     Let $d\geq 2$, $0<s<1$ and $A\in L^2_{\rm loc}(\R^d)$. Then there exists
     a function $k_A$ on $\R^d\times\R^d$ such that the following
     holds. If $\chi_0,\ldots,\chi_n$ are Lipschitz continuous functions
     on $\R^d$ satisfying $\sum_{j=0}^n \chi_j^2 \equiv 1$, then one has
     \begin{equation}\label{eq:imsmag}
       \| |D-A|^s u\|^2
       = \sum_{j=0}^n \||D-A|^s \chi_j u\|^2 - (u,L_A u),
       \qquad u\in\dom |D-A|^s,
     \end{equation}
     where $L_A$ is the bounded operator with integral kernel
     \begin{equation*}
       L_A(x,y) := k_A(x,y) \sum_{j=0}^n (\chi_j(x)-\chi_j(y))^2.
     \end{equation*}
     Moreover, for a.e. $x,y\in\R^d$
     \begin{equation*}
       |k_A(x,y)| \leq a_{s,d} |x-y|^{-d-2s},
       \qquad {\rm and\ hence\ }\quad
       |L_A(x,y)| \leq L(x,y)
     \end{equation*}
with $L$ defined in Lemma~\ref{ims}.
\end{lemma}

\begin{proof}
     By the argument of \cite{Si1} we can choose a form core for
     $|D-A|^{2s}$ which is invariant under multiplication by Lipschitz
     continuous functions. It suffices to prove \eqref{eq:imsmag} only
     for functions $u$ from such a core.

     We write $k_A(x,y,t) := \exp(-t|D-A|^{2s})(x,y)$ for the heat kernel
     and find
     \begin{align*}
       & \sum_{j=0}^n  (\chi_j u, (1-\exp(-t|D-A|^{2s})) \chi_j u)
       = (u, (1-\exp(-t|D-A|^{2s})) u)\\
       & \qquad\qquad  + \frac 12\sum_{j=0}^n
       \iint k_A(x,y,t)(\chi_j(x)-\chi_j(y))^2 \overline{u(x)} u(y)
       \,dx\,dy.
     \end{align*}
     Now we divide by $t$ and note that by our assumption on $u$ the
     left side
     converges to $\sum_{j=0}^n \||D-A|^s \chi_j u\|^2$ as $t\to
     0$. Similarly the first term on the right side divided by $t$ converges to
     $\||D-A|^s u\|^2$. Hence the last term divided by $t$ converges to
     some limit $(u, L_A u)$. The diamagnetic inequality \eqref{eq:diamagheat}
     yields the bound $|k_A(x,y,t)| \leq
     \exp(-t(-\Delta)^{2s})(x,y)$. This implies in particular that $L_A$
     is a bounded operator. Now it is easy to check that $L_A$ is an
     integral operator and that the absolute value of its kernel is
     bounded pointwise by the one of $L$.
\end{proof}

The following helps to clarify the role of the kernel $k_A$.

\begin{corollary}\label{locop}
     Let $u\in\dom|D-A|^{2s}$ and assume that $\Omega:=\R^d\setminus\supp
     u \not=\emptyset$. Then
     \begin{equation*}
       \left(|D-A|^{2s} u\right)(x) = - \int_{\R^d} k_A(x,y) u(y)\, dy
       \qquad \textrm{for\ $x\in\Omega$}.
     \end{equation*}
\end{corollary}

\begin{proof}
     Let $\phi\in C_0^\infty(\Omega)$ and choose $\chi_0, \chi_1$ such
     that $\chi_0\equiv 1$ on $\supp u$, $\chi_1\equiv 1$ on $\supp\phi$
     and $\chi_0^2+\chi_1^2\equiv 1$. By polarization, \eqref{eq:imsmag}
     implies $(\phi,|D-A|^{2s}u) = -(\phi,L_A u)= -\int \overline{\phi(x)} 
k_A(x,y)
     u(y) \,dx\,dy$, whence the assertion.
\end{proof}


\section{Closability of the Quadratic Form
     $b_{\beta,A}$}\label{app:closable}

This appendix contains some technical details concerning the quadratic
form $b_{\beta,A}$ defined in (\ref{def:bba}). In particular, we shall
show its closability.
Throughout this appendix we assume that $d\geq 2$, $0<s\leq 1$ and
$A\in L^2_{\rm loc}(\R^d)$.

\begin{lemma}\label{formcore}
     The sets $C_0^\infty(\R^d\setminus\{0\})$ and $\mathcal D := \{ w\in
     \dom(D-A)^2 \cap L^\infty(\R^d) : \ \supp w \text{ compact in }
     \R^d\setminus\{0\} \}$ are form cores for $|D-A|^{2s}$.
\end{lemma}

\begin{proof}
     It suffices to prove the statement for $s=1$. In this case it is
     proved in \cite{Si1} that $C_0^\infty(\R^d)$ and $\mathcal D^* := \{
     w\in \dom(D-A)^2 \cap L^\infty(\R^d) : \ \supp w \text{ compact} \}$
     are form cores for $(D-A)^2$. Hence the statement will follow if we
     can approximate every function in any of these two spaces by
     functions from the same space vanishing in a neighborhood of the
     origin. But for functions $u$ from $C_0^\infty(\R^d)$ or $\mathcal
     D^*$ both functions $Du$ and $Au$ are square-integrable. This
     reduces the lemma to the case $A=0$ where it is well-known.
\end{proof}

Now let $d-2s<\beta<d$ and recall that the quadratic form $q_{\beta,A}$
was defined in (\ref{def:qba}). Note that
$C_0^\infty(\R^d\setminus\{0\})$ is invariant under the unitary
transformation in (\ref{def:qbu}). Therefore, closability of $q_{\beta,A}$
and $b_{\beta,A}$ on $C_0^\infty(\R^d\setminus\{0\})$ are equivalent.

\begin{lemma}\label{closable}
     The quadratic form $q_{\beta,A}$, defined in (\ref{def:qba}), is closable 
on
     $C_0^\infty(\R^d\setminus\{0\})$.
\end{lemma}

\begin{proof}
     It suffices to show closability of the form 
$r_{\beta,A}[u]:=h_{s,A}[|x|^\alpha u]$
     on $C_0^\infty(\R^d\setminus\{0\})$
     for $0< \alpha = (\beta+2s-d)/2 < s$.

     Let $\mathcal D$ be as in Lemma \ref{formcore}. We shall show
     that the quadratic form $r_{\beta,A}$ on $\mathcal D$ is closable
     and that $C_0^\infty(\R^d\setminus\{0\})$ is dense with respect to
     $(r_{\beta,A}[w]+\|w\|^2)^{1/2}$ in the closure of $\mathcal D$
     with respect to this norm. This implies the assertion.

     Let $w\in\mathcal D$. Since $|x|^\alpha$ is smooth on $\supp w$ and
     $\dom(D-A)^2$ is invariant under multiplication by smooth functions
     we have $|x|^\alpha w\in\dom(D-A)^2$ and hence $|x|^\alpha
     w\in\dom|D-A|^{2s}$. Since $|x|^\alpha w$ has compact support it
     follows from Lemma \ref{imsmag} and Corollary \ref{locop} that
     $|(|D-A|^{2s}|\cdot|^\alpha w)(x)|  \leq C(w) |x|^{-d-2s}$ for all large 
$|x|$. In particular,
     $|x|^\alpha |D-A|^{2s}|x|^\alpha w \in L^2(\R^d)$. Moreover,
     $|x|^{-2(s-\alpha)}w\in L^2(\R^d)$. It follows that if
     $u_n\in\mathcal D$ such that $u_n\to 0$ in $L^2(\R^d)$, then the
     bilinear form associated with $r_{\beta,A}$ satisfies
     \begin{equation*}
       r_{\beta,A}[u_n,v]
       = (u_n, |x|^\alpha |D-A|^{2s}|x|^\alpha w
       - \mathcal C_{s,d}|x|^{-2(s-\alpha)}w ) \to 0
     \end{equation*}
     as $n\to\infty$. By standard arguments, this proves that
     $r_{\beta,A}$ is closable on $\mathcal D$.

     In order to show the density of $C_0^\infty(\R^d\setminus\{0\})$ we
     let again $w\in\mathcal D$. Since $|x|^\alpha w\in\dom|D-A|^s$,
     Lemma \ref{formcore} yields a sequence $u_n\in
     C_0^\infty(\R^d\setminus\{0\})$ such that $\| |D-A|^s
     (u_n-|x|^\alpha w)\| + \|u_n-|x|^\alpha w\| \to 0$. Hence, if $w_n
     := |x|^{-\alpha} u_n$, then $w_n\to w$ in $L^2(\R^d\setminus\mathcal
     B)$ and $0 \leq r_{\alpha,A}[w_n-w] \leq \| |D-A|^s |x|^\alpha
     (w_n-w)\|^2 \to 0$. Moreover, Hardy's inequality implies also that $\|
     |x|^{-s+\alpha} (w_n-w)\| \to 0$ and hence $w_n\to w$ in
     $L^2(\mathcal B)$. This proves the density of
     $C_0^\infty(\R^d\setminus\{0\})$ in the closure of $\mathcal D$ with
     respect to $(r_{\alpha,A}[w]+\|w\|^2)^{1/2}$.
\end{proof}


\section{Nash's Argument}\label{app:nash}

For the sake of completeness we recall here how the Nash inequality
\eqref{eq:weightednash} and the contraction property of $\exp(-t B_\beta)$
imply the heat kernel estimate \eqref{eq:heat}. Let $k(t):= \|\exp(-t
B_\beta) v\|_{\mathfrak H_\beta}^2$. Then
\begin{equation*}
\begin{split}
k'(t) & =  -2 b_\beta[\exp(-t B_\beta) v] \\
& \leq - 2 C_{q,d,s}'^{-1} k(t)^{1+1/p}
\|\exp(-tB_\beta) v\|_{L^1(|x|^{-\beta}dx)}^{-2/p} \\
& \leq - 2 C_{q,d,s}'^{-1} k(t)^{1+1/p}
\|v\|_{L^1(|x|^{-\beta}dx)}^{-2/p},
\end{split}
\end{equation*}
where we used \eqref{eq:weightednash} and the contraction property. Hence
\begin{equation*}
\left(k(t)^{-1/p}\right)' \geq 2 p^{-1} C_{q,d,s}'^{-1}
\|v\|_{L^1(|x|^{-\beta}dx)}^{-2/p}
\end{equation*}
and, after integration,
\begin{equation*}
k(t)^{-1/p} \geq 2 p^{-1} C_{q,d,s}'^{-1}
t \|v\|_{L^1(|x|^{-\beta}dx)}^{-2/p}.
\end{equation*}
This means that
\begin{equation*}
\|\exp(-t B_\beta) v\|_{\mathfrak H_\beta}^2
\leq (p C_{q,d,s}'/ 2)^{p} t^{-p} \|v\|_{L^1(|x|^{-\beta}dx)}^2.
\end{equation*}
By duality, noting that $B_\beta$ is self-adjoint in
$L^2(|x|^{-\beta}dx)$, this also implies
\begin{equation*}
\|\exp(-t B_\beta) v\|_{L^\infty(|x|^{-\beta}dx)}^2
\leq (p C_{q,d,s}'/ 2)^{p} t^{-p} \|v\|_{\mathfrak H_\beta}^2.
\end{equation*}
Finally, by the semi-group property $\exp(-t B_\beta)=\exp(-t
B_\beta/2)\exp(-t B_\beta/2)$
\begin{equation*}
\|\exp(-t B_\beta) v\|_{L^\infty(|x|^{-\beta}dx)}^2
\leq (p C_{q,d,s}'/ 2)^{2p} (t/2)^{-{2p}} \|v\|_{L^1(|x|^{-\beta}dx)}^2.
\end{equation*}
This is exactly the estimate \eqref{eq:heat} with the constant given in
Proposition \ref{heat}.


\section{A Trace Estimate}\label{app:trace}

In this appendix, we sketch the argument leading to
Proposition~\ref{trace}. We emphasize that we shall ignore several
technical details. We assume that $W$ is smooth with compact support
in $\R^d\setminus\{0\}$, and we put $k(x,y,t):=\exp(-t
B_{\beta})(x,y)$. We claim that the trace formula
\begin{equation}\label{eq:traceformula}
     \begin{split}
       & \tr_{\mathfrak H_\beta} F(W^{1/2}B_{\beta}^{-1}W^{1/2})
       \\ & = \int_0^\infty \frac{dt}t  \lim_{n\to\infty}
       \int_{\R^d}\cdots\int_{\R^d}
       \prod_{j=1}^n k\left(x_j,x_{j-1},\frac tn\right)
       f\left(\frac tn \sum_{j=1}^n W(x_j)\right)
       \, \frac{dx_1}{|x_1|^{\beta}} \cdots \frac{dx_n}{|x_n|^{\beta}}
     \end{split}
\end{equation}
holds true for any non-negative, lower semi-continuous function $f$
vanishing near the origin. Here, $F$ is related to $f$ as in (\ref{eq:f}). In 
the
integral in \eqref{eq:traceformula} we use the convention that
$x_0=x_n$. Indeed, by an approximation argument it suffices to prove
this formula for
\begin{equation*}
     F(\lambda)=\lambda/(1+\alpha\lambda)
     \quad,\quad
     f(\mu)=\mu e^{-\alpha \mu}\,,
\end{equation*}
where $\alpha>0$ is a constant. Using the resolvent identity and
Trotter's product formula, one easily verifies that
\begin{equation*}
     \begin{split}
       F(W^{1/2}B_\beta^{-1}W^{1/2})
       & = W^{1/2}(B_\beta+\alpha W)^{-1}W^{1/2} \\
       & = \int_0^\infty
       W^{1/2}\exp(-t(B_\beta+\alpha W)) W^{1/2}\,dt \\
       & = \int_0^\infty \lim_{n\to\infty} T_n(t) \,dt
     \end{split}
\end{equation*}
in this case. Here,
\begin{equation*}
     T_n(t) :=
     W^{1/2}
     \big(\exp(-tB_\beta/n)\exp(-t\alpha W/n)\big)^n
     W^{1/2}.
\end{equation*}
The latter is an integral operator. We evaluate its trace via
integrating its kernel on the diagonal. Then we arrive at
\begin{equation*}
     \tr_{\mathfrak H_\beta} T_n(t) =
      \idotsint
       \prod_{j=1}^n k\left(x_j,x_{j-1},\frac tn\right)
       W(x_n) e^{-\frac{\alpha t}n \sum_j W(x_j)}
       \, \frac{dx_1}{|x_1|^{\beta}} \cdots \frac{dx_n}{|x_n|^{\beta}}\,.
\end{equation*}
After symmetrization with respect to the variables this leads to
\begin{equation*}
     \tr_{\mathfrak H_\beta} T_n(t) = \frac 1t
     \idotsint
     \prod_{j=1}^n k\left(x_j,x_{j-1},\frac tn\right)
     f\left(\frac tn \sum_j W(x_j)\right)
     \, \frac{dx_1}{|x_1|^{\beta}} \cdots \frac{dx_n}{|x_n|^{\beta}}\,.
\end{equation*}
The claimed formula \eqref{eq:traceformula} follows if one
interchanges the trace with the $t$-integration and the
$n$-limit.

Now we assume in addition that $f$ is convex. Then Jensen's inequality yields
\begin{equation*}
     \begin{split}
       & \idotsint
       \prod_{j=1}^n k\left(x_j,x_{j-1},\frac tn\right)
       f\left(\frac tn \sum_j W(x_j)\right) \,  \frac{dx_1}{|x_1|^{\beta}} 
\cdots \frac{dx_n}{|x_n|^{\beta}} \\
       & \qquad \leq \idotsint
       \prod_{j=1}^n k\left(x_j,x_{j-1},\frac tn\right)
       \frac 1n \sum_j f(tW(x_j)) \,  \frac{dx_1}{|x_1|^{\beta}} \cdots 
\frac{dx_n}{|x_n|^{\beta}} \\
       & \qquad = \idotsint
       \prod_{j=1}^n k\left(x_j,x_{j-1},\frac tn\right)
       f(tW(x_1)) \,  \frac{dx_1}{|x_1|^{\beta}} \cdots 
\frac{dx_n}{|x_n|^{\beta}}\,.
     \end{split}
\end{equation*}
(Eq. (\ref{eq:traceformula}) holds also in the magnetic case discussed
in Section~\ref{sec:diamag}. Before applying
Jensen's inequality, one first has to use the diamagnetic inequality 
(\ref{eq:ba}) to eliminate
the magnetic field in the kernel $k$, however.)
Finally, we use the semigroup property to integrate with respect to
the variables $x_2,\ldots,x_n$. We find that the latter integral is
equal to
\begin{equation*}
     \int_{\R^d} k(x,x,t) f(tW(x)) \, \frac{dx}{|x|^\beta}\,.
\end{equation*}
Plugging this into \eqref{eq:traceformula} leads to the estimate
\eqref{eq:trace}.

Details concerning the justification of the above manipulations can be
found in \cite{RoSo} (see Theorem 2.4 there).

\end{appendix}


\bibliographystyle{amsalpha}

\end{document}